\documentclass[11pt]{amsart}

\usepackage{amssymb,bbold,mathtools}
\usepackage{mathrsfs}
\usepackage{mathtools}
\usepackage{enumitem}
\usepackage{hyperref}
\usepackage{cleveref}
\usepackage{verbatim}

\theoremstyle{plain}
\newtheorem{theorem}{Theorem}[section]
\newtheorem{proposition}[theorem]{Proposition}
\newtheorem{corollary}[theorem]{Corollary}
\newtheorem{lemma}[theorem]{Lemma}

\theoremstyle{definition}
\newtheorem{remark}[theorem]{Remark}
\newtheorem{example}[theorem]{Example}
\newtheorem{definition}[theorem]{Definition}

\newcommand{\abs}[1]{\lvert#1\rvert}
\newcommand{\norm}[1]{\lVert#1\rVert}

\newcommand{\uppars}[1]{\textup{(}#1\textup{)}}

\newcommand{\RR}{\mathbb R}

\newcommand{\pos}[1]{{#1}^+}
\newcommand{\posE}{\pos{E}}
\newcommand{\posF}{\pos{F}}
\newcommand{\nega}[1]{{#1}^-}

\newcommand{\ob}{\mathrm{ob}}
\newcommand{\oc}{\mathrm{oc}}
\newcommand{\reg}{\mathrm r}
\newcommand{\obops}{\mathcal{L}_{\ob}}
\newcommand{\ocops}{\mathcal{L}_{\oc}}
\newcommand{\regops}{\mathcal{L}_{\reg}}

\newcommand{\odual}[1]{{#1}^{\thicksim}}
\newcommand{\ocdual}[1]{{#1}^\thicksim_{\oc}}

\newcommand{\Ell}{\mathrm{L}}
\newcommand{\Cee}{\mathrm{C}}

\newcommand{\conv}{\xrightarrow}
\newcommand{\convwithoverset}[1]{\conv{#1}}

\renewcommand{\o}{\mathrm{o}}

\newcommand{\uo}{\mathrm{uo}}

\newcommand{\unb}{\mathrm{u}}

\newcommand{\uoLt}{{\widehat{\tau}}}

\newcommand{\oconv}{\convwithoverset{\o}}
\newcommand{\uoconv}{\convwithoverset{\uo}}

\newcommand{\tauconv}{\convwithoverset{\tau}}

\newcommand{\uoLtconv}{\convwithoverset{\uoLt}}

\newcommand{\netgen}[2]{{(#1)_{#2}}}

\newcommand{\is}{{\mathcal A}}
\newcommand{\net}{\netgen{x_\alpha}{\alpha\in\is}}

\newcommand{\istwo}{{\mathcal B}}

\newcommand{\seq}[2]{(#1)_{#2=1}^\infty}
\newcommand{\seqxn}{\seq{x}{n}}

\newcommand{\oLtop}{o-Lebes\-gue topology}
\newcommand{\oLtops}{o-Lebes\-gue topologies}

\newcommand{\uoLtop}{uo-Lebes\-gue topology}
\newcommand{\uoLtops}{uo-Lebes\-gue topologies}

\newcommand{\oadhtext}{o-ad\-her\-ence}
\newcommand{\uoadhtext}{uo-ad\-her\-ence}
\newcommand{\soadhtext}{$\sigma$o-ad\-her\-ence}
\newcommand{\suoadhtext}{$\sigma$uo-ad\-her\-ence}

\newcommand{\adh}[2]{a_{#1}(#2)}
\newcommand{\sadh}[2]{a_{\sigma #1}(#2)}

\newcommand{\oadh}[1]{\adh{\o}{#1}}
\newcommand{\uoadh}[1]{\adh{\uo}{#1}}
\newcommand{\soadh}[1]{\sadh{\o}{#1}}
\newcommand{\suoadh}[1]{\sadh{\uo}{#1}}

\newcommand{\oclos}[1]{\overline{#1}^{\o}}
\newcommand{\uoclos}[1]{\overline{#1}^{\uo}}

\newcommand{\suoclos}[1]{\overline{#1}^{\sigma\uo}}

\setlist[enumerate,1]{label=\textup{(\arabic*)},ref=\textup{(\arabic*)}}
\setlist[enumerate,2]{label=\textup{(\alph*)},ref=\textup{(\alph*)}}
\setlist[enumerate,3]{label=\textup{(\roman*)},ref=\textup{(\roman*)}}
\setlist[enumerate,4]{label=\textup{(\Alph*)},ref=\textup{(\Alph*)}}

\newlist{enumerate_arabic}{enumerate}{1}
\setlist[enumerate_arabic,1]{label=\textup{\arabic*.},ref=\textup{\arabic*}}

\newlist{enumerate_alpha}{enumerate}{1}
\setlist[enumerate_alpha,1]{label=\textup{(\alph*)},ref=\textup{(\alph*)}}

\newlist{enumerate_roman}{enumerate}{1}
\setlist[enumerate_roman,1]{label=\textup{(\roman*)},ref=\textup{(\roman*)}}

\newlist{enumerate_Alpha}{enumerate}{1}
\setlist[enumerate_Alpha,1]{label=\textup{(\Alph*)},ref=\textup{(\Alph*)}}

\crefname{theorem}{Theorem}{Theorems}
\crefname{proposition}{Proposition}{Propositions}
\crefname{lemma}{Lemma}{Lemmas}
\crefname{corollary}{Corollary}{Corollaries}
\crefname{conjecture}{Conjecture}{Conjectures}
\crefname{definition}{Definition}{Definitions}
\crefname{example}{Example}{Examples}
\crefname{remark}{Remark}{Remarks}
\crefname{assumption}{Assumption}{Assumptions}
\crefname{hypothesis}{Hypothesis}{Hypotheses}
\crefname{question}{Question}{Questions}
\crefname{problem}{Problem}{Problems}
\crefname{task}{Task}{Tasks}
\crefname{addendum}{Addendum}{Addenda}
\crefname{idea}{Idea}{Ideas}
\crefname{suggestion}{Suggestion}{Suggestions}
\crefname{context}{Context}{Contexts}
\crefname{exercise}{Exercise}{Exercises}
\crefname{section}{Section}{Section}
\crefname{subsection}{Section}{Section}
\crefname{subsubsection}{Section}{Section}

\crefname{equation}{equation}{equations}

\usepackage[charter]{mathdesign}

\numberwithin{theorem}{section}
\numberwithin{equation}{section}

\allowdisplaybreaks 

\begin{document}
\title[Vector lattices with a Hausdorff uo-Lebesgue topology]{Vector lattices with a Hausdorff uo-Lebesgue topology}

\author{Yang Deng}
\address{Yang Deng; School of Economic Mathematics, Southwestern University of Finance and Economics, Chengdu, Sichuan 611130, People's Republic of China}
\email{dengyang@swufe.edu.cn}

\author{Marcel de Jeu}
\address{Marcel de Jeu; Mathematical Institute, Leiden University, P.O.\ Box 9512, 2300 RA Leiden, the Netherlands;
	and Department of Mathematics and Applied Mathematics, University of Pretoria, Cor\-ner of Lynnwood Road and Roper Street, Hatfield 0083, Pretoria, South Africa}
\email{mdejeu@math.leidenuniv.nl}


\begin{abstract}
	We investigate the construction of a Hausdorff uo-Lebesgue topology on a vector lattice from a Hausdorff (o)-Lebesgue topology on an order dense ideal, and what the properties of the topologies thus obtained are. When the vector lattice has an order dense ideal with a separating order continuous dual, it is always possible to supply it with such a topology in this fashion, and the restriction of this topology to a regular sublattice is then also a Hausdorff uo-Lebesgue topology. A regular vector sublattice of $\mathrm{L}_0(X,\Sigma,\mu)$ for a semi-finite measure $\mu$ falls into this category, and the convergence of nets in its Hausdorff uo-Lebesgue topology is then the convergence in measure on subsets of finite measure. When a vector lattice not only has an order dense ideal with a separating order continuous dual, but also has the countable sup property, we show that every net in a regular vector sublattice that converges in its Hausdorff uo-Lebesgue topology always contains a sequence that is uo-convergent to the same limit. This enables us to give satisfactory answers to various topological questions about uo-convergence in this context.
\end{abstract}

\maketitle

\section{Introduction and overview}

\noindent In this paper, we investigate the construction of a Hausdorff uo-Lebesgue topology on a vector lattice from a Hausdorff (o)-Lebesgue topology\footnote{In the literature, what we call a \oLtop\ is simply called a Lebesgue topology. Now that \uoLtops,  with a completely analogous definition, have become objects of a more extensive study, it seems consistent to also add a prefix to the original term.}\ on an order dense ideal, and what the properties of the topologies thus obtained are.

After recalling the relevant notions and making the necessary preparations in \cref{sec:preliminaries}, the key construction is carried out in  \cref{res:generated_topology} in \cref{sec:unbounded_topologie_generated_by_topologies_on_ideals}, below. The idea of starting with a topology on an order dense ideal originates from \cite{conradie:2005} but, whereas the construction in \cite{conradie:2005} to obtain a global topology is carried out using Riesz pseudo-norms, we follow an approach using neighbourhood bases of zero that is inspired by \cite{taylor_THESIS:2018, taylor:2019}. Using such neighbourhood bases, it is possible to perform the construction under minimal hypotheses on the initial data, and thus understand how these hypotheses are reflected in the properties of the resulting global topology. The remainder of \cref{sec:unbounded_topologie_generated_by_topologies_on_ideals} is mainly concerned with showing how the general theorem relates to existing results in the literature. Our working with neighbourhood bases of zero enables us to explain certain `pathologies' in the literature, where a topology of unbounded type is not Hausdorff, or not linear, from the general theorem.

In \cref{sec:uoLtops_going_up_and_going_down}, we move to the context where the initial ideal is actually order dense and admits a Hausdorff \oLtop. In that case, every regular vector sublattice of the global vector lattice admits a Hausdorff \uoLtop. The resulting overview \cref{res:overview}, below, mostly consists of a summary of results that are already in the literature, though not presented in this way. It is also recalled in that section that a regular vector sublattice admits a Hausdorff \uoLtop\ when the global vector lattice admits one. Consequently, there is a going-up-going-down procedure: starting with a Hausdorff \oLtop\ on an order dense ideal, one obtains a Hausdorff \uoLtop\ on the global vector lattice, and then finally also one on every regular vector sublattice.

In view of the going-up-going-down construction, it is evidently desirable to have a class of vector lattices that admit Hausdorff \oLtops\ because such data can serve as `germs' for Hausdorff \uoLtops. The vector lattices with  separating order continuous duals form such a class, and this is exploited in \cref{sec:separating_order_continuous_dual}.

\cref{sec:vector_lattices_of_equivalence_classes_of_measurable_functions} is concerned with regular vector sublattices of $\Ell_0(X,\Sigma,\mu)$ for a semi-finite measure $\mu$. Via the going-up-going-down principle, every regular vector sublattice of  $\Ell_0(X,\Sigma,\mu)$ admits a Hausdorff \uoLtop. We give a rigorous proof of the fact that the convergence of nets in such a topology is the convergence in measure on subsets of finite measure. For $\Ell_p(X,\Sigma,\mu)$, we also discuss how the (in fact) unique Hausdorff \uoLtop\ on these spaces can be described in various seemingly different ways that are still equivalent. The relation between these topologies and minimal and smallest Hausdorff locally solid linear topologies on these spaces is explained.

\cref{sec:uo-convergent_subsequences_of_uoLt-convergent_nets} is concerned with convergent sequences that can always be found `within' nets that are convergent in a Hausdorff \uoLtop\ on a vector lattice that has the countable sup property and that has an order dense ideal with a separating order continuous dual. The precise statement is in \cref{res:tau_m_to_sub_uo}, below; this is one of the main theorems in this paper. It is in the same spirit as the fact that a sequence that converges (globally) in measure always contains a subsequence that converges almost everywhere to the same limit.

Finally, in \cref{sec:topological aspects of uo-convergence}, we study topological aspects of uo-convergence. The relations between uo-convergence and various order topologies are not at all well understood, but when the global vector lattice has the countable sup property, and also has an order dense ideal with a separating order continuous dual, then a reasonably satisfactory picture emerges. In \cref{res:seven_sets_equal} and \cref{res:eleven_sets_equal}, below, various topological closures and (sequential) adherences are then seen to be equal. It is then also possible to give a necessary and sufficient criterion for sequential uo-convergence to be topological; see \cref{res:sequential_uo_convergence_topological_two}, below.

\smallskip

We have tried to be as complete in the development of this part of the theory of uo-convergence as we could, and also to relate to relevant existing results in the literature whenever possible. Any omissions at this point are unintentional.

\section{Preliminaries}\label{sec:preliminaries}

\noindent In this section, we collect a number of definitions, notations, conventions and preparatory results. We refer the reader to the textbooks \cite{abramovich_aliprantis_INVITATION_TO_OPERATOR_THEORY:2002}, \cite{aliprantis_border_INFINITE_DIMENSIONAL_ANALYSIS_THIRD_EDITION:2006}, \cite{aliprantis_burkinshaw_LOCALLY_SOLID_RIESZ_SPACES_WITH_APPLICATIONS_TO_ECONOMICS_SECOND_EDITION:2003},  \cite{aliprantis_burkinshaw_POSITIVE_OPERATORS_SPRINGER_REPRINT:2006},  \cite{luxemburg_zaanen_RIESZ_SPACES_VOLUME_I:1971}, \cite{meyer-nieberg_BANACH_LATTICES:1991}, \cite{schaefer_BANACH_LATTICES_AND_POSITIVE_OPERATORS:1974}, \cite{zaanen_RIESZ_SPACES_VOLUME_II:1983}, and \cite{zaanen_INTRODUCTION_TO_OPERATOR_THEORY_IN_RIESZ_SPACES:1997} for general background information on vector lattices and Banach lattices.

\subsection{Vector lattices, operators, and (unbounded) order convergence}

\noindent All vector spaces are over the real numbers. Measures take their values in $[0,\infty]$ and are not supposed to satisfy any condition unless otherwise specified. All vector lattices are supposed to be Archimedean. The positive cone of a vector lattice $E$ is denoted by $\posE$.

Let $E$ be a vector lattice, and let $F$ be a vector sublattice of $E$. Then $F$ is \emph{order dense in $E$} when, for every $x\in E$ with $x>0$, there exists a $y\in F$ such that $0<y\leq x$; $F$ is called \emph{super order dense in $E$} when, for every $x\in \posE$, there exists a sequence $\seqxn\subseteq \posF$ with $x_n\uparrow x$ in $E$. The vector sublattice $F$ of $E$ is order dense in $E$ if and only if, for every $x\in \posE$, we have $x=\sup\{y\in F: 0\leq y\leq x\}$; see \cite[Theorem~1.34]{aliprantis_burkinshaw_POSITIVE_OPERATORS_SPRINGER_REPRINT:2006}, for example.

A vector sublattice $F$ of a vector lattice $E$ is called \emph{majorising in $E$} when, for every $x\in E$, there exists a $y\in F$ such that $x\leq y$. In some sources, such as \cite{conradie:2005}, $F$ is then said to be full in $E$.

A vector lattice $E$ \emph{has the countable sup property} when, for every non-empty subset $S$ of $E$ that has a supremum in $E$, there exists an at most countable subset of $S$ that has the same supremum in $E$ as $S$. In parts of the literature, such as in \cite{luxemburg_zaanen_RIESZ_SPACES_VOLUME_I:1971} and \cite{zaanen_INTRODUCTION_TO_OPERATOR_THEORY_IN_RIESZ_SPACES:1997}, $E$ is then said to be order separable. A vector lattice $E$ has the countable sup property if and only if, whenever a net $\net\subseteq\pos{E}$ and $x\in\pos{E}$ are such that $x_\alpha\uparrow x$ in $E$, there exists a sequence of indices $\seq{\alpha_n}{n}$ in $\is$ such that $x_{\alpha_n}\uparrow x$ in $E$; see \cite[Theorem~23.2.(iii)]{luxemburg_zaanen_RIESZ_SPACES_VOLUME_I:1971}.

Let $E$ be a vector lattice, and let $x\in E$. We say that a net $\net$ in $E$ is \emph{order convergent to $x\in E$} (denoted by $x_\alpha\oconv x$) when there exists a net $(y_\beta)_{\beta\in \mathscr{B}}$ in $E$ such that $y_\beta\downarrow 0$ and with the property that, for every $\beta_0\in\istwo$, there exists an $\alpha_0\in \is$ such that $\abs{x-x_\alpha}\leq y_{\beta_0}$ whenever $\alpha$ in $\is$ is such that $\alpha\geq\alpha_0$.  Note that the index sets $\is$ and $\istwo$ need not be equal; for a discussion of the difference between these two possible definitions we refer to \cite{abramovich_sirotkin:2005}, for example.

Let $E$ and $F$ be vector lattices. The order bounded operators from $E$ into $F$ will be denoted by $\obops(E,F)$, and the regular operators from $E$ into $F$ by $\regops(E,F)$. When $F$ is Dedekind complete, we have $\obops(E,F)=\regops(E,F)$, and this space is then a Dedekind complete vector lattice; see \cite[Theorem~1.18]{aliprantis_burkinshaw_POSITIVE_OPERATORS_SPRINGER_REPRINT:2006}, for example. We write $\odual{E}$ for $\obops(E,\RR)=\regops(E,\RR)$.

A linear operator  $T: E\to F$ between two vector lattices $E$ and $F$ is \emph{order continuous} when, for every net $\net$ in $E$,  the fact that $x_\alpha\oconv 0$ in $E$ implies that $Tx_\alpha\oconv 0$ in $F$. When $T$ is positive one can, equivalently, require that, for every net $\net$ in $E$, the fact that $x_\alpha\downarrow 0$ in $E$ imply that $Tx_\alpha\downarrow 0$ in $F$. An order continuous linear operator between two vector lattices is automatically order bounded; see \cite[Lemma~1.54]{aliprantis_burkinshaw_POSITIVE_OPERATORS_SPRINGER_REPRINT:2006}, for example. The order continuous linear operators from $E$ into $F$ will be denoted by $\ocops(E,F)$. In the literature, the notation $\mathcal L_{\mathrm n}(E,F)$ is often used. When $F$ is Dedekind complete, $\ocops(E,F)$ is a band in $\regops(E,F)$; see \cite[Theorem~1.57]{aliprantis_burkinshaw_POSITIVE_OPERATORS_SPRINGER_REPRINT:2006}, for example. We write $\ocdual{E}$ for $\ocops(E,\RR)$.

The following result is easily established using the Riesz-Kantorovich formulas and their `dual versions'; see \cite[Theorems~1.18 and~1.23]{aliprantis_burkinshaw_POSITIVE_OPERATORS_SPRINGER_REPRINT:2006}, for example. We shall be interested only in the case where the lattice $F$ in it is the real numbers and the band $B$ is the zero band, but the general case comes at no extra cost in the routine proof.

\begin{proposition}\label{res:polars}
	Let $E$ and $F$ be vector lattices, where $F$ is Dedekind complete, and let $B$ be a band in $F$.
	\begin{enumerate}
		\item Let $I$ be an ideal of $E$. Then the subset
		\[
		\{T\in\regops(E,F): Tx \in B \text{ for all }x\in I\}
		\]
		of $\regops(E,F)$ is band in $\regops(E,F)$. For every subset $S$ of $I$ that generates $I$, it is equal to
		\[
		\{T\in\regops(E,F): \abs{T}\abs{x} \in B \text{ for all }x\in S\}.
		\]
		\item Let $\mathcal I$ be an ideal of $\regops(E,F)$. Then the subset
		\[
		\{x\in E: Tx \in B \text{ for all }T\in \mathcal I\}
		\]
		of $E$ is an ideal of $E$. For every subset $\mathcal S$ of $\mathcal I$ that generates $\mathcal I$, it is equal to
		\[
		\{x\in E: \abs{T}\abs{x} \in B \text{ for all }T\in \mathcal S\}.
		\]
		It is a band in $E$ when $\mathcal I\subseteq\ocops(E,F)$.
	\end{enumerate}	
\end{proposition}

Let $F$ be a vector sublattice of a vector lattice $E$. Then $F$ is a \emph{regular vector sublattice of $E$} when the inclusion map from $F$ into $E$ is order continuous.
Equivalently, for every net $\net$ in $F$, the fact that $x_\alpha\downarrow 0$ in $F  $ should imply that $x_\alpha\downarrow 0$ in $E$. It is immediate from the latter criterion that ideals are regular vector sublattices. It is also true that order dense vector sublattices are regular vector sublattices; see \cite[Theorem~1.23]{aliprantis_burkinshaw_LOCALLY_SOLID_RIESZ_SPACES_WITH_APPLICATIONS_TO_ECONOMICS_SECOND_EDITION:2003}, for example.

Let $\net$ be a net in a vector lattice $E$, and let $x\in E$. We say that $(x_\alpha)$ is \emph{unbounded order convergent to $x$ in $E$}  (denoted by $x_\alpha\uoconv x$) when $\abs{x_\alpha-x}\wedge y\oconv 0$ in $E$ for all $y\in \posE$. Order convergence implies unbounded order convergence to the same limit. For order bounded nets, the two notions coincide. \footnote{Although we shall not need this, it would be less than satisfactory not to mention here that the uo-continuous dual of a vector lattice (defined in the obvious way) has a very concrete description, and is often trivial. According to \cite[Proposition~2.2]{gao_leung_xanthos:2018}, it is the linear span of the coordinate functionals corresponding to atoms.}

We shall repeatedly refer to the following collection of results; see \cite[Theorem~2.8, Corollary~2.12, and Theorem~3.2]{gao_troitsky_xanthos:2017}.

\begin{theorem}\label{res:local-global_for_o-convergence_and_uo-convergence}
	Let $E$ be a vector lattice, and let $F$ be a vector sublattice of $E$. Take a net $\net$ in $F$.
	
	\begin{enumerate}
		\item Suppose that $F$ is order dense and majorising in $E$. Then $x_\alpha\oconv 0$ in $F$ if and only if $x_\alpha\oconv 0$ in $E$.
		\item\label{part:order_convergence_and_regular_sublattices} Suppose that $F$ is a regular vector sublattice of $E$ and that $\net$ is order bounded in $F$. Then $x_\alpha\oconv 0$ in $F$ if and only if $x_\alpha\oconv 0$ in $E$.
		\item\label{part:unbounded_order_convergence_and_regular_sublattices} The following are equivalent:
		\begin{enumerate}
			\item $F$ is a regular vector sublattice of $E$;
			\item for every net $\net$ in $F$, the fact that $x_\alpha\uoconv 0 $ in $F$ implies that $x_\alpha\uoconv 0$ in $E$;
			\item for every net $\net$ in $F$, $x_\alpha\uoconv 0 $ in $F$ if and only if $x_\alpha\uoconv 0$ in $E$.
		\end{enumerate}	
	\end{enumerate}		
\end{theorem}

In the sequel of this paper, we shall encounter restrictions of order continuous linear functionals on vector lattices to vector sublattices. For this, we include the following result. It is based on a theorem of Veksler's. It contains quite a bit more than we shall actually need, but we use the opportunity to present the results in it, and its fourth and fifth parts in particular.

\begin{theorem}\label{res:veksler}
	Let $E$ be a vector lattice, let $F$ be a vector sublattice of $E$, and let $G$ be a Dedekind complete vector lattice. Take $T\in \ocops(E,G)$.
	\begin{enumerate}
		\item Suppose that $F$ is a regular vector sublattice of $E$. Then the restriction $T|_F:F\to G$ of $T$ to $F$ is order continuous.\label{part:veksler_1}
		\item Suppose that $F$ is a regular sublattice of $E$. When $\ocops(E,G)$ separates the points of $E$, then $\ocops(F,G)$ separates the points of $F$.\label{part:veksler_2}
		\item Suppose that $F$ is an order dense vector sublattice of $E$.  Then the restriction map $T\mapsto T|_F$ is a positive linear injection from $\ocops(E,G)$ into $\ocops(F,G)$.\label{part:veksler_3}
	\end{enumerate}
		
	Suppose that $F$ is an order dense and majorising vector sublattice of $E$. Then:
	\begin{enumerate}[resume]
		\item the restriction map $T\mapsto T|_F$ is a lattice isomorphism between $\ocops(E,G)$ and $\ocops(F,G)$;\label{part:veksler_4}
		\item $\ocops(E,G)$ separates the points of $E$ if and only if $\ocops(F,G)$ separates the points of $F$.\label{part:veksler_5}
	\end{enumerate}
\end{theorem}

\begin{proof}
	Part~\ref{part:veksler_1} is obvious, and then so is part~\ref{part:veksler_2}.
	
	It is clear from part~\ref{part:veksler_1}  that $\ocops(F,G)$ separates the points of $F$ whenever $\ocops(E,G)$ separates the points of $E$.
	
	Suppose that $F$ is an order dense (hence regular) vector sublattice of $E$ and that $T\in\ocops(E,G)$ is such that $T|_F=0$. Take $x\in \posE$. Then $\{y\in F: 0\leq y\leq x\}\uparrow x$ in $E$. Since $T|_F=0$, the order continuity of $T$ on $E$ then implies that $Tx=0$. Hence $T=0$, and we conclude that the restriction map $T\mapsto T|_F$ is a positive linear injection from $\ocops(E,G)$ into $\ocops(F,G)$.

	Suppose that $F$ is order dense and majorising in $E$.
	
	Take $S\in\ocops(F,G)$. In that case, according to a result of Veksler's (see \cite[Theorem~1.65]{aliprantis_burkinshaw_POSITIVE_OPERATORS_SPRINGER_REPRINT:2006}), each of $\pos{S}$ and $\nega{S}$ can be extended to a positive order continuous operator from $E$ into $G$. Hence $S$ itself can be extended to an order continuous operator $S^{\text{ext}}$ from $E$ into $G$. By what we have already observed in part~\ref{part:veksler_3}, such an order continuous extension is unique, and we conclude from this that the map $S\mapsto S^{\text{ext}}$ is a positive linear injection  from $\ocops(F,G)$ into $\ocops(E,G)$. It is clear that the extension and restriction maps between $\ocops(E,G)$ and $\ocops(F,G)$ are each other's inverses. We conclude that the restriction map $T\mapsto T|_F$ is a bi-positive linear bijection between $\ocops(E,G)$ and $\ocops(F,G)$. Hence it is a lattice isomorphism, as required.
	
	One direction of the equivalence in part~\ref{part:veksler_5} is clear from part~\ref{part:veksler_2}. For the converse direction, suppose that $\ocops(F,G)$ separates the points of $F$. Take $x\in E$ such that $Tx=0$ for all $T\in \ocops(E,G)$. Since $\ocops(E,G)$ is an ideal of $\regops(E,F)$, \cref{res:polars} shows that $T\abs{x}=0$ for all $T\in \ocops(E,G)$. Suppose that $x\neq 0$. Then there exists a $y\in F$ such that $0<y\leq\abs{x}$, and we have  $Ty=0$ for all positive $T\in \ocops(E,G)$, hence for all $T\in \ocops(E,G)$. In view of part~\ref{part:veksler_4}, this is the same as saying that $Sy=0$ for all $S\in\ocops(F,G)$. Our assumption yields that $y=0$; this contradiction shows that we must have $x=0$.
\end{proof}

\subsection{Topologies on vector lattices}

\noindent When $E$ is a vector space, a \emph{linear topology on $E$} is a (not necessarily Hausdorff) topology that provides $E$ with the structure of a topological vector space. When $E$ is a vector lattice, a \emph{locally solid linear topology on $E$} is a linear topology on $E$ such that there exists a base of (not necessarily open) neighbourhoods of 0 that are solid subsets of $E$.
For the general theory of locally solid linear topologies on vector lattices we refer to \cite{aliprantis_burkinshaw_LOCALLY_SOLID_RIESZ_SPACES_WITH_APPLICATIONS_TO_ECONOMICS_SECOND_EDITION:2003}. A locally solid linear topology on $E$ that is also a locally convex linear topology is a \emph{locally convex-solid linear topology}.
In that case, there exists a base of neighbourhoods of 0 that consists of absorbing, closed, convex, and solid subsets of $E$; see \cite[p.~59]{aliprantis_burkinshaw_LOCALLY_SOLID_RIESZ_SPACES_WITH_APPLICATIONS_TO_ECONOMICS_SECOND_EDITION:2003}.

When $E$ is a vector lattice, a \emph{locally solid additive topology on $E$} is a topology that provides the additive group $E$ with the structure of a (not necessarily Hausdorff) topological group, such that there exists a base of (not necessarily open) neighbourhoods of 0 that are solid subsets of $E$.

Let $E$ be a vector lattice. We say that \emph{order convergence in $E$ is topological} when there exists a (evidently unique) topology on $E$ such that its convergent nets are precisely the order convergent nets, with preservation of limits. It follows from the properties of order convergence that such a topology is automatically a Hausdorff linear topology. Likewise, \emph{unbounded order convergence in $E$ is topological} when there exists a topology on $E$ such that its convergent nets are precisely the nets that are unbounded order convergent, with preservation of limits. Such a topology is again unique, and automatically a Hausdorff linear topology.

A topology $\tau$ on a vector lattice $E$ is an \emph{\oLtop} when it is a (not necessarily Hausdorff) locally solid linear topology on $E$ such that, for a net $\net$ in $E$ and $x\in E$, the fact that $x_\alpha\oconv x$ in $E$ implies that $x_\alpha\tauconv x$. Equivalently, the fact that $x_\alpha\oconv 0$ in $E$ should imply that $x_\alpha\tauconv 0$.
A vector lattice need not admit a Hausdorff \oLtop. It can be shown, see \cite[Example~3.2]{aliprantis_burkinshaw_LOCALLY_SOLID_RIESZ_SPACES_WITH_APPLICATIONS_TO_ECONOMICS_SECOND_EDITION:2003}, that $\Cee([0,1])$ does not even admit a Hausdorff locally solid linear topology such that \emph{sequential} order convergence implies topological convergence.

A topology $\tau$ on a vector lattice $E$ is a \emph{\uoLtop} when it is a (not necessarily Hausdorff) locally solid linear topology on $E$ such that, for a net $\net$ in $E$ and $x\in E$, the fact that $x_\alpha\uoconv x$ in $E$ implies that $x_\alpha\tauconv x$. Equivalently, the fact that $x_\alpha\uoconv 0$ in $E$ should imply that $x_\alpha\tauconv 0$.
Since order convergence implies unbounded order convergence, a \uoLtop\ is an \oLtop.

The following fundamental facts are from \cite[Proposition~3.2,~3.4, and~6.2, and Corollary~6.3]{conradie:2005} and
\cite[Theorems~5.5,~5.9, and~6.4]{taylor:2019}.

\begin{theorem}[Conradie and Taylor]\label{res:conradie_taylor}
	Let $E$ be a vector lattice. Then the following are equivalent:
	\begin{enumerate}
		\item $E$ admits a Hausdorff \oLtop;
		\item $E$ admits a Hausdorff \uoLtop;\label{part:uoLtop}
		\item the partially ordered set of all Hausdorff locally solid linear topologies on $E$ has a minimal element.\label{part:mintop}
	\end{enumerate}
When this is the case, the topologies in the parts~\ref{part:uoLtop} and~\ref{part:mintop} are both unique, they coincide, and they are the smallest Hausdorff \oLtop\ on $E$.
\end{theorem}

When $E$ admits a Hausdorff \uoLtop, we shall denote the unique such topology by $\uoLt_E$. In \cite{conradie:2005}, it is denoted by $\tau_m$. For a given vector lattice, there may be several ways to obtain a Hausdorff \uoLtop\ on it. This can then give criteria for the convergence of nets in the common resulting topology that are apparently equivalent, but not always immediately obviously so. See \cref{rem:various_descriptions} for this, for example.

\begin{remark}\label{rem:minimal_and_smallest}
	Some caution is necessary when consulting the literature on minimal Hausdorff locally solid linear topologies because in \cite[Definition~7.64]{aliprantis_burkinshaw_LOCALLY_SOLID_RIESZ_SPACES_WITH_APPLICATIONS_TO_ECONOMICS_SECOND_EDITION:2003} such a topology is defined as what would usually be called a \emph{smallest} Hausdorff locally solid linear topologies. When a vector lattice $E$ admits a complete metrisable \oLtop, such as a Banach lattice with an order continuous norm, then it admits a smallest (in the usual sense of the word) Hausdorff locally solid linear topology; see \cite[Theorem~7.65]{aliprantis_burkinshaw_LOCALLY_SOLID_RIESZ_SPACES_WITH_APPLICATIONS_TO_ECONOMICS_SECOND_EDITION:2003}. Combining this with \cref{res:conradie_taylor}, we see that $E$ then admits a (necessarily unique) Hausdorff \uoLtop\ $\uoLt_E$, and that $\uoLt_E$ is then not just the smallest Hausdorff \oLtop, but even the smallest Hausdorff locally solid linear topology on $E$.
\end{remark}

\section{Unbounded topologies generated by topologies on ideals}\label{sec:unbounded_topologie_generated_by_topologies_on_ideals}

\noindent We shall now describe how topologies `of unbounded type' on vector lattices can be obtained from topologies on ideals. There are already several constructions in this vein and accompanying results in the literature; see \cite{conradie:2005,deng_o_brien_troitsky:2017,kandic_li_troitsky:2018,kandic_marabeh_troitsky:2017,taylor_THESIS:2018,taylor:2019}, for example. In the following result, we carry out such a construction in what appears to be the most general possible context. Starting from a locally solid (not necessarily linear or Hausdorff) additive topology on an ideal $F$ of a vector lattice $E$ and a non-empty subset of $F$, we define an `unbounded' locally solid additive topology on $E$. Various known results in more special cases can then be understood from the general theorem, as will be discussed in \crefrange{exam:taylor:2019}{exam:conradie:2005}, below.

It is important to note that, in our construction, the initial topology on the ideal $F$ need not be the restriction of a global topology on $E$. For such restricted topologies, several results are already available in \cite[Propositions~9.3 and~9.4]{taylor:2019}.  It will, however, become clear in \cref{rem:essential_difference}, below, that such a global topology (of a suitable type) need not always exist. Hence there is an actual gain by starting from topologies on ideals, and a concrete illustration of this procedure can be found in  \cref{sec:vector_lattices_of_equivalence_classes_of_measurable_functions}. The possibility of working in such a general setting was already observed in \cite[Remark~2.8, Example~2.9, and Remark~2.10]{taylor_THESIS:2018}, where it was also noted that the continuity of the scalar multiplication may fail to hold in the new topology. Part~\ref{part:last} of \cref{res:generated_topology}, below, gives necessary and sufficient conditions for this continuity.

The subset $S$ figuring in the construction can be replaced by the ideal that it generates without altering the result. Although it may conceptually be more natural to work with ideals than with subsets, working with arbitrary subsets has the advantage of keeping an eye on a small number of presumably relatively easily manageable `test elements'. It is for this reason that we carry this along to later results; see also \cref{rem:small_subset_criterion}, below. The convenience of this approach will become apparent in the proof of \cref{res:tau_m_is_convergence_in_measure}.

\begin{theorem}\label{res:generated_topology}
	Let $E$ be a vector lattice, let $F$ be an ideal of $E$, and let $\tau_F$ be a \uppars{not necessarily Hausdorff} locally solid additive topology on $F$. Take a non-empty subset $S$ of $F$.
	
	There exists a unique \uppars{possibly non-Hausdorff} additive topology $\unb_S\!\tau_{F}$ on $E$ such that, for a net $\net$ in $E$, $x_\alpha\conv{\unb_S\!\tau_{F}} 0$ in $E$ if and only if $\abs{x_\alpha}\wedge \abs{s}\conv {\tau_F}0$ in $F$ for all $s\in S$.
	
	Let $I_S\subseteq F$ be the ideal generated by $S$ in $E$. For a net $\net$ in $E$,	$x_\alpha\conv{\unb_S\!\tau_{F}} 0$ in $E$ if and only if $\abs{x_\alpha}\wedge \abs{y}\conv {\tau_F}0$ in $F$ for all $y\in I_S$.
	
	Furthermore:
	\begin{enumerate}
		\item\label{part:first} the inclusion map from $F$ into $E$ is $\tau_F$--$\unb_S\!\tau_F$-continuous;
		\item the topology $\unb_S\!\tau_{F}$ on $E$ is a locally solid additive topology;
		\item the following are equivalent:
		\begin{enumerate}
			\item $\unb_S\!\tau_{F}$ is a Hausdorff topology on $E$;\label{test}
			\item $\tau_F$ is a Hausdorff topology on $F$ and $I_S$ is order dense in $E$;
		\end{enumerate}
	\item\label{part:last} the following are equivalent:
	\begin{enumerate}
		\item\label{part:linear_condition_i} for all $x\in E$ and $s\in S$,
		\begin{equation}\label{eq:absorbing_condition}
		\abs{\varepsilon x}\wedge \abs{s}\conv{\tau_F} 0
		\end{equation}
		in $F$ as $\varepsilon\to 0$ in $\mathbb R$;
		\item\label{part:linear_condition_ii} for all $x\in E$ and $y\in I_S$, $\abs{\varepsilon x}\wedge \abs{y}\conv{\tau_F} 0$ in $F$
		as $\varepsilon\to 0$ in $\mathbb R$;
		\item\label{part:linear_condition_iii} $\unb_S\!\tau_F$ is a \uppars{possibly non-Hausdorff} linear topology on $E$.
	\end{enumerate}
	
	\end{enumerate}
\end{theorem}

\begin{proof}
	Suppose that $\tau_F$ is a (not necessarily Hausdorff) locally solid additive topology on $F$.
	
	The uniqueness of $\unb_S\!\tau_{F}$ is clear because the nets converging to 0 and then, by translation invariance of the topology, to arbitrary points of $E$ are prescribed.
	
	We turn to the existence of such a topology $\unb_S\!\tau_{F}$. Take a neighbourhood base $\{U_\lambda\}_{\lambda\in\Lambda}$ of zero in $F$ for $\tau_F$ consisting of solid subsets of $F$. For $y\in I_S$ and $\lambda\in\Lambda$, set
	\begin{equation}\label{eq:local_base_1}
	V_{\lambda,y}\coloneqq\{x\in E: \abs{x}\wedge\abs{y}\in U_\lambda\}.
	\end{equation}
	The $V_{\lambda,y}$ are solid subsets of $E$ since $F$ is an ideal of $E$ and the $U_\lambda$ are solid subsets of $F$. Set
	\begin{equation}\label{eq:local_base_2}
	\mathcal{N}_0\coloneqq \{V_{\lambda,y}: \lambda\in \Lambda, y\in I_S \}.
	\end{equation}
	 We claim that $\mathcal{N}_0$ is a base of neighbourhoods of zero for a topology on $E$, which we shall already denote by $\unb_S\!\tau_{F}$, that provides the additive group $E$ with the structure of a topological group.  Necessary and sufficient conditions on $\mathcal N_0$ for this can be found in \cite[Theorem~3 on p.46]{husain_INTRODUCTION-TO_TOPOLOGICAL_GROUPS:1966}; we now verify these.
	
	 Take $V_{\lambda_1,y_1},V_{\lambda_2,y_2}\in \mathcal{N}_0$. There exists a $\lambda_3\in\Lambda$ such that $U_{\lambda_3}\subseteq U_{\lambda_1}\cap U_{\lambda_2}$. Take $x\in V_{\lambda_3,\abs{y_1}\vee\abs{y_2}}$. Then
	 \[
	 \abs{x}\wedge\abs{y_1}\leq\abs{x}\wedge\left({\abs{y_1}\vee\abs{y_2}}\right)\in U_{\lambda_3}\subseteq U_{\lambda_1}.
	 \]
	 Since $F$ is an ideal of $E$ and $U_{\lambda_1}$ is a solid subset of $F$, this implies that $\abs{x}\wedge\abs{y_1}\in U_{\lambda_1}$, so that $x\in V_{\lambda_1,y_1}$. Likewise, $x\in V_{\lambda_2,y_2}$, and we see that $V_{\lambda_3,\abs{y_1}\vee\abs{y_2}}\subseteq V_{\lambda_1,y_1}\cap V_{\lambda_2,y_2}$.
	
	 It is evident that $V_{\lambda,y}=-V_{\lambda,y}$ for all $V_{\lambda,y}\in\mathcal N_0$.
	
	 Take $V_{\lambda,y}\in\mathcal N_0$. There exists a $\mu\in\Lambda$ such that $ U_\mu + U_\mu \subseteq U_\lambda$. Then, for all $x_1,x_2\in V_{\mu,y}$, we have
	\[
	\abs{x_1+x_2}\wedge \abs{y}\leq \abs{x_1}\wedge\abs{y}+\abs{x_2}\wedge\abs{y}\in U_\mu+U_\mu\subseteq U_\lambda.
	\]
	Since $F$ is an ideal of $E$ and $U_\lambda$ is a solid subset of $F$, this implies that $\abs{x_1+x_2}\wedge \abs{y}\in U_{\lambda}$, so that $x_1+x_2\in V_{\lambda,y}$. Hence $V_{\mu,y}+V_{\mu,y}\subseteq V_{\lambda,y}$.
	
	An appeal to \cite[p.~46, Theorem~3]{husain_INTRODUCTION-TO_TOPOLOGICAL_GROUPS:1966} now establishes our claim.
	
It is clear from the definition of $\unb_S\!\tau_F$ that, for a net $\net$ in $E$, $x_\alpha\conv{\unb_S\!\tau_{F}} 0$ in $E$ if and only if $\abs{x_\alpha}\wedge \abs{y}\conv {\tau_F}0$ in $F$ for all $y\in I_S$; this observation goes back to \cite[Remark~2.3]{taylor_THESIS:2018}.
	
	Certainly, the fact that $\abs{x_\alpha}\wedge \abs{y}\conv {\tau_F}0$ in $F$ for all $y\in I_S$ implies that  $\abs{x_\alpha}\wedge \abs{s}\conv {\tau_F}0$ in $F$ for all $s\in S$. Conversely, suppose that $\net$ is a net in $E$ such that $\abs{x_\alpha}\wedge \abs{s}\conv {\tau_F}0$ in $F$ for all $s\in S$. Take $y\in I_S$. There exist $s_1,\dotsc,s_n\in S$ and integers $k_1,\dotsc,k_n\geq 1$ such that $\abs{y}\leq \sum_{i=1}^{n} k_i\abs{s_i}$. Hence $\abs{x_\alpha}\wedge \abs{y}\leq \sum_{i=1}^{n} k_i \left(\abs{x_\alpha}\wedge \abs{s_i}\right)$.  Since $\tau_F$ is a locally solid additive topology on $F$, this implies that $\abs{x_\alpha}\wedge \abs{y}\conv{\tau_F} 0$ in $F$.
	
	We turn to the parts~\ref{part:first}--\ref{part:last}.
	
	Since $F$ is an ideal of $E$ and the $U_\lambda$ are solid subsets of $F$, we have $U_\lambda\subseteq V_{\lambda,y}$ for all $\lambda\in\Lambda$ and $y\in I_S$. This implies that the inclusion map from $F$ into $E$ is $\tau_F$--$\unb_S\!\tau_F$-continuous.
	
	The topology $\unb_S\!\tau_F$ is a locally solid additive topology on $E$ by construction.
	
	Suppose that $\unb_S\!\tau_{F}$ is a Hausdorff topology on $E$. Then so is the topology it induces on $F$, which is weaker than $\tau_F$. Hence $\tau_F$ is a Hausdorff topology on $F$. Take $x\in E$ with $x>0$. Then there exists a $V_{\lambda,y}\in \mathcal{N}_0$ with $x\notin V_{\lambda,y}$. In particular, $x\wedge\abs{y}\neq 0$. Hence $0<x\wedge\abs{y} \leq x$. Since  $x\wedge\abs{y}\in I_S$, we see that $I_S$ is order dense in $E$.
	
	Suppose, conversely, that $\tau_F$ is a Hausdorff topology on $F$ and that $I_S$ is order dense in $E$. Take $x\neq 0$ in $E$. There exists a $y\in I_S$ with $0<y\leq \abs{x}$. Pick $U_{\lambda_0}\in \{U_\lambda\}_{\lambda\in \Lambda}$ such that $y\notin U_{\lambda_0}$. Then $\abs{x}\wedge \abs{y}=y\notin U_{\lambda_0}$, so that $x\notin V_{\lambda_0,y}$. Hence $\bigcap_{V\in\mathcal N_0}V=\{0\}$. By \cite[p.~48, Theorem~4]{husain_INTRODUCTION-TO_TOPOLOGICAL_GROUPS:1966}, $\unb_S\!\tau_F$ is a Hausdorff additive topology on the topological group $E$.
	
	We shall now verify the equivalence of the parts~\ref{part:linear_condition_i}--\ref{part:linear_condition_iii} of part~\ref{part:last}.
	
	We prove that part~\ref{part:linear_condition_i} implies part~\ref{part:linear_condition_ii}. Take $x\in E$ and $y\in I_S$. There exist $s_1,\dotsc,s_n\in S$ and integers $k_1,\dotsc,k_n\geq 1$ such that $\abs{y}\leq \sum_{i=1}^{n} k_i\abs{s_i}$, and it follows from this that $\abs{\varepsilon x}\wedge \abs{y}\leq \sum_{i=1}^{n}  k_i\left(\abs{\varepsilon x}\wedge \abs{s_i}\right)$ for all $\varepsilon\in\RR$.  Since $\tau_F$ is a locally solid additive topology on $F$, it follows that $\abs{\varepsilon x}\wedge \abs{y}\conv{\tau_F} 0$ in $F$ as $\varepsilon\to 0$ in $\RR$.
	
We prove that part~\ref{part:linear_condition_ii} implies part~\ref{part:linear_condition_iii}. Fix $\lambda\in\Lambda$ and $y\in I_S$, and take $x\in E$. Since $\abs{\varepsilon x}\wedge \abs{y}\conv{\tau_F} 0$ in $F$ as $\varepsilon\to 0$ in $\RR$, there exists a $\delta>0$ such that $\abs{\varepsilon x}\wedge \abs{y}\in U_\lambda$ whenever $\abs{\varepsilon}\leq\delta$. That is, $\varepsilon x\in V_{\lambda,y}$ whenever $\abs{\varepsilon}\leq\delta$. This implies that $V_{\lambda,y}$ is absorbing. Furthermore, since $V_{\lambda,y}$ is a solid subset of $E$, it is clear that $\varepsilon x\in V_{\lambda,y}$ whenever $x\in V_{\lambda,y}$ and $\varepsilon\in[-1,1]$. Hence $V_{\lambda,y}$ is balanced.
	Then \cite[Theorem 5.6]{aliprantis_border_INFINITE_DIMENSIONAL_ANALYSIS_THIRD_EDITION:2006} implies that $\unb_S\!\tau_F$ is a linear topology on $E$.
	
We prove that part~\ref{part:linear_condition_iii} implies part~\ref{part:linear_condition_i}. Take $x\in E$. Then $\varepsilon x\conv{\unb_S\!\tau_F} 0$ in $E$ as $\varepsilon\to 0$ in $\RR$. By construction, this implies (and is, in fact, equivalent to) the fact that  $\abs{\varepsilon x}\wedge \abs{s}\conv {\tau_F}0$ in $F$ for all $s\in S$.

This concludes the proof of the equivalence of the three parts of part~\ref{part:last}. The proof of the theorem is now complete.
	
\end{proof}

\begin{definition}
	The topology $\unb_S\!\tau_{F}$ in \cref{res:generated_topology,res:generated_topology} is called the \emph{unbounded topology on $E$ that is generated by $\tau_F$ via $S$}.
\end{definition}

\begin{remark}\label{rem:two_sets_give_the_same_topology}
	It is clear from the two equivalent criteria  in \cref{res:generated_topology} for a net in $E$ to be $\unb_S\!\tau_F$-convergent to zero that $\unb_S\!\tau_F=\unb_{I_S}\!\tau_F$ for every non-empty subset $S$ of $F$. Consequently, $\unb_{S_1}\!\tau_F=\unb_{S_2}\!\tau_F$ whenever $S_1,S_2$ are non-empty subsets of $F$ such that $I_{S_1}=I_{S_2}$.
\end{remark}

\begin{remark}
	In \cref{res:generated_topology}, suppose that the locally solid additive topology $F$ is the restriction $\tau_E|_F$ of a locally solid additive topology on $E$. It is then easy to see that $\unb_S\!\left(\tau_E|_F\right)=\unb_S\!\tau_E$ for every non-empty subset $S$ of $F$.
\end{remark}

\begin{remark}\label{rem:explicit_neighbourhood_base}
	In \cref{res:generated_topology}, and also in the remainder of this paper, the topologies of interest are characterised by their convergent nets. It should be noted, however, that in \cref{eq:local_base_1,eq:local_base_2} the proof of \cref{res:generated_topology} provides an explicit neighbourhood base of zero in $E$ for $\unb_S\!\tau_F$, in terms of a neighbourhood base of zero in $F$ for $\tau_F$ and the ideal $I_S$. Suppose, for example that $\tau_F$ is a (possibly non-Hausdorff) locally convex-solid linear topology on $F$ that is generated by a family $\{\rho_\gamma: \gamma\in \Gamma\}$ of lattice semi-norms on $F$, as will be the case in \cref{sec:separating_order_continuous_dual}. Then the collection of subsets of $E$ of the form
	\[
	\{x\in E: \rho_i(\abs{x}\wedge\abs{y})<\varepsilon\text{ for } \rho_1,\dotsc\rho_n\in\Gamma\},
	\]
	where $y\in I_S$, $n\geq 1$, and $\varepsilon>0$ are arbitrary, is a neighbourhood base of zero in $E$ for $\unb_S\!\tau_F$.
\end{remark}

Our next result is concerned with iterating the construction in \cref{res:generated_topology}. It generalises what is in \cite[p.~997]{taylor:2019}.

\begin{proposition}\label{res:repeating_the_construction}
	Let $E$ be a vector lattice, let $F_1$ be an ideal of $E$, and let $\tau_{F_1}$ be a \uppars{not necessarily Hausdorff} locally solid additive topology on $F_1$. Take a non-empty subset $S_1$ of $F_1$, and consider the unbounded topology $\unb_{S_1}\!\!\tau_{F_1}$ on $E$ that is generated by $\tau_{F_1}$ via $S_1$. Let $F_2$ be an ideal of $E$, and let $\left(\unb_{S_1}\!\!\tau_{F_1}\right)|_{F_2}$ denote the topology on $F_2$ that is induced on $F_2$ by $\unb_{S_1}\!\!\tau_{F_1}$. Then $\left(\unb_{S_1}\!\!\tau_{F_1}\right)|_{F_2}$ is a locally solid additive topology on $F_2$. Take a non-empty subset $S_2$ of $F_2$. Then $\unb_{S_2}\!\left[\left(\unb_{S_1}\!\!\tau_{F_1}\right)|_{F_2}\right]=\unb_{I_{S_1}\cap I_{S_2}}\!\!\tau_{F_1}$.
	In particular, when $S$ is a non-empty subset of $F_1\cap F_2$, then  $\unb_{S}\!\left[\left(\unb_{S}\!\!\tau_{F_1}\right)|_{F_2}\right]=\unb_{S}\!\!\tau_{F_1}$.
\end{proposition}

\begin{proof}
It is clear from \cref{res:generated_topology} that $\left(\unb_{S_1}\tau_{F_1}\right)|_{F_2}$ is a locally solid additive topology on $F_2$. Let $\net$ be a net in $E$. Then we have the following chain of equivalent statements:
\begin{align*}
&x_\alpha\conv{\unb_{S_2}\left[\left(\unb_{S_1}\tau_{F_1}\right)|_{F_2}\right]}0 \text{ in }E\\\
&\iff \abs{x_\alpha}\wedge \abs{y_2}\conv{\left(\unb_{S_1}\tau_{F_1}\right)|_{F_2}}0 \text{ in }F_2\text{ for all }y_2\in I_{S_2}\\
&\iff \abs{x_\alpha}\wedge \abs{y_2}\conv{\unb_{S_1}\tau_{F_1}   }0 \text{ in }E\text{ for all }y_2\in I_{S_2}\\
&\iff \abs{x_\alpha}\wedge \abs{y_2}\wedge\abs{y_1}\conv{\tau_{F_1}}0 \text{ in }F_1\text{ for all }y_1\in I_{S_1}\text{ and } y_2\in I_{S_2}\\
&\iff \abs{x_\alpha}\wedge \abs{y}\conv{\tau_{F_1}}0 \text{ in }F_1\text{ for all }y\in I_{S_1}\cap I_{S_2}\\
&\iff x_\alpha\conv{\unb_{I_{S_1}\cap I_{S_2}}\tau_{F_1}}0 \text{ in }E.
\end{align*}
Hence $\unb_{S_2}\!\left[\left(\unb_{S_1}\!\!\tau_{F_1}\right)|_{F_2}\right]=\unb_{I_{S_1}\cap I_{S_2}}\!\!\tau_{F_1}$.
\end{proof}

\begin{remark}
	In \cref{res:repeating_the_construction}, suppose that $\tau_{F_1}$ is a \uppars{not necessarily Hausdorff} locally solid additive topology on $F_1$ such that, for all $x\in E$ and $s\in S_1$, $\abs{\varepsilon x}\wedge \abs{s}\conv{\tau_{F_1}} 0$ in $F_1$ as $\varepsilon\to 0$ in $\mathbb R$. It is then clear from \cref{res:generated_topology} that $\unb_{S_1} \!\!\tau_{F_1}$, $\left(\unb_{S_1}\!\!\tau_{F_1}\right)|_{F_2}$, and $\unb_{I_{S_1}\cap I_{S_2}}\!\!\tau_{F_1}$  are \uppars{possibly non-Hausdorff} locally solid linear topologies on $E$, $F_2$, and $E$, respectively.
	
\end{remark}

We shall now explain how \cref{res:generated_topology} relates to various results already in the literature.

\begin{example}\label{exam:taylor:2019}
	 When $F=E$ and $\tau_E$ is a locally solid linear topology on $F=E$, the condition in \cref{eq:absorbing_condition} is automatically satisfied for any non-empty subset $S$ of $F=E$. According to \cref{res:generated_topology}, $\unb_E\!\tau_E$ is a locally solid linear topology on $E$ that is Hausdorff if and only if $\tau_E$ is Hausdorff; this is \cite[Theorem~2.3]{taylor:2019}. Furthermore, when $A$ is an ideal of $E$, $\unb_A\!\tau_E$ is a locally solid linear topology on $E$ that is Hausdorff if and only if $\tau_E$ is Hausdorff and $A$ is order dense in $E$; this is \cite[Propositions~9.3 and ~9.4]{taylor:2019}.
\end{example}

\begin{example}\label{exam:deng_o'brien_troitsky:2018}
	Let $E$ be a Banach lattice. In \cref{res:generated_topology}, we take $F=E$, for $\tau_F$ we take the norm topology $\tau_E$ on $F=E$, and for $S\subseteq F$ we take $S=F=E$. Then the condition in \cref{eq:absorbing_condition} is satisfied. According to \cref{res:generated_topology}, $\unb_E\!\tau_E$ is a Hausdorff locally solid linear topology on $E$ and, for a net $\net$ in $E$, $x_\alpha\conv{\unb_E\!\tau_E}0$ if and only if $\norm{\abs{x_\alpha}\wedge\abs{y}}\to 0$ for all $y\in E$. In \cite{deng_o_brien_troitsky:2017}, this type of convergence is called \emph{unbounded norm convergence}, or \emph{un-convergence} for short. It was already observed in \cite[Section~7]{deng_o_brien_troitsky:2017} that it is topological; in \cite[p.~746]{kandic_li_troitsky:2018},  $\unb_F\!\tau_F$ is then called the \emph{un-topology}.
\end{example}
	
\begin{example}\label{exam:kandic_li_troitsky:2018}
	Let $E$ be a vector lattice, and let $F$ be an ideal of $E$ that is a normed vector lattice. In \cref{res:generated_topology}, we take for $\tau_F$ the norm topology on $F$, and for $S\subseteq F$ we take $S=F$. According to \cref{res:generated_topology}, $\unb_F\!\tau_F$ is a (possibly non-Hausdorff) additive topology on $E$ and, for a net $\net$ in $E$, $x_\alpha\conv{\unb_F\!\tau_F}0$ if and only if $\norm{\abs{x_\alpha}\wedge\abs{y}}\to 0$ for all $y\in F$. This type of convergence is called \emph{un-convergence with respect to $X$} in \cite{kandic_li_troitsky:2018}. It was already observed that it is topological in \cite[p.~747]{kandic_li_troitsky:2018}, where $\unb_F\!\tau_F$ is called the \emph{un-topology on $E$ induced by $F$}.
	
	In \cite[Example~1.3]{kandic_li_troitsky:2018}, it is shown that  $\unb_F\!\tau_F$ can fail to be a Hausdorff topology on $E$. Since $\tau_F$ is a Hausdorff topology on $F$, \cref{res:generated_topology} shows that the pertinent ideal $F$ in \cite[Example~1.3]{kandic_li_troitsky:2018} must fail to be order dense in $E$; this is indeed easily seen to be the case. \cref{res:generated_topology} implies that $\unb_F\!\tau_F$ is Hausdorff if and only if $F$ is order dense in $F$; this is \cite[Proposition~1.4]{kandic_li_troitsky:2018}.
	
	In \cite[Example~1.5]{kandic_li_troitsky:2018}, it is shown that $\unb_F\!\tau_F$ can fail to be a linear topology on $E$. According to \cref{res:generated_topology}, the condition in \cref{eq:absorbing_condition} must fail to be satisfied in the context of \cite[Example~1.5]{kandic_li_troitsky:2018}; this is indeed easily seen to be the case. \cref{res:generated_topology} shows that $\unb_F\!\tau_F$ always provides $E$ with an additive topology; this was also noted in \cite[p.~748]{kandic_li_troitsky:2018} in that particular context.
	
	In \cite[p.~748]{kandic_li_troitsky:2018}, the authors observe that $\unb_F\!\tau_F$ is a locally solid linear topology on the vector lattice $E$ whenever $E$ is a normed lattice and the norm on $E$ extends that on $F$, and also whenever the norm on $F$ is order continuous. Both facts follow from \cref{res:generated_topology} because \cref{eq:absorbing_condition} is then satisfied. This is clear when $E$ is a normed lattice and the norm on $E$ extends that on $F$. Suppose that the norm on $F$ is order continuous. Take $x\in E$ and $y\in F$. Then $\abs{\varepsilon x}\wedge \abs{y}\oconv 0$ in $E$ as $\varepsilon\to 0$. Since the net $ \abs{\varepsilon x}\wedge \abs{y}$ is order bounded in the ideal $F$ of $E$, which is a regular vector sublattice of $E$, \cref{res:local-global_for_o-convergence_and_uo-convergence} implies that $\abs{\varepsilon x}\wedge \abs{y}\oconv 0$ in $F$, and then $\abs{\varepsilon x}\wedge \abs{y}\conv{\tau_F}0$ as $\varepsilon\to 0$ by the order continuity of the norm on $F$.
\end{example}

\begin{example}\label{exam:zabeti:2018}
For a vector lattice $E$, we let $\abs{\sigma}(E,\odual{E})$ denote its absolute weak topology; the definition of this locally solid linear topology will be recalled in \cref{sec:separating_order_continuous_dual}. Taking $E=F=S$ in \cref{res:generated_topology} yields the so-called \emph{unbounded absolute weak topology $\unb_E\abs{\sigma}(E,\odual{E})$} on $E$. It is a locally solid additive topology on $E$ that is Hausdorff if and only if $\odual{E}$ separates the points of $E$. When $E$ is a Banach lattice, $\unb_E\abs{\sigma}(E,\odual{E})$ is a Hausdorff locally solid linear topology on $E$. It is studied in \cite{zabeti:2018}.	
	
\end{example}

\begin{example}\label{exam:conradie:2005}
	In \cite[p.~290]{conradie:2005}, a construction is given to obtain a locally solid linear topology on a vector lattice $E$ from a locally solid linear topology on an ideal $F$ of $E$. This is done using Riesz pseudo-norms, rather than by working with neighbourhood bases of zero as we have done. The key ingredient is to start with a Riesz pseudo-norm $p$ on $F$, take an element $u$ of $\posF$, and introduce a map $p_u:E\to\RR$ by setting $p_u(x)\coloneqq p(\abs{x}\wedge u)$ for $x\in E$.  It is then remarked that $p_u$ is a Riesz pseudo-norm on $E$. This need not always be the case, however. By way of counter-example, take for $E$ the vector lattice of all real-valued functions on $\RR$, and for $F$ the ideal of $E$ consisting of all bounded functions on $\RR$. For $p$, we take the supremum norm on $F$. For $u\in\posF$, we choose the constant function 1. We define $x\in E$ by setting $x(t)\coloneqq t$ for $t\in\RR$. Then $p_u(\lambda x)=\norm{\abs{\lambda x}\wedge u}=1$ for all non-zero $\lambda\in\RR$, whereas we should have that $\lim_{\lambda\to 0} p_u(\lambda x)=0$. This implies that the topologies on $E$ that are thus constructed, although locally solid additive topologies, need not be linear topologies. This `pathology' is similar to that in \cite[Example~1.5]{kandic_li_troitsky:2018} that was mentioned above; our example here is also quite similar to that in \cite[Example~1.5]{kandic_li_troitsky:2018}. Fortunately, in the continuation of the argument in \cite{conradie:2005}, $p$ is taken to be a Riesz pseudo-norm on $F$ that is continuous with respect to a Hausdorff \oLtop\ $\tau_F$ on $F$. In this context, $p_u$ \emph{is} a Riesz pseudo-norm on $E$. Indeed, since $F$, being an ideal of $E$, is a regular vector sublattice of $E$, \cref{res:local-global_for_o-convergence_and_uo-convergence} easily yields that $\abs{\lambda x}\wedge u\oconv  0$ in $F$ as $\lambda\to 0$. Since $\tau_F$ is an \oLtop\ on $E$, we have $\abs{\lambda x}\wedge u\conv{\tau_F} 0$ in $F$ as $\lambda\to 0$, and then the continuity of $p$ on $F$ yields that $p_u(\lambda x)\to 0$ as $\lambda\to 0$. Thus the construction in \cite{conradie:2005} proceeds correctly after all. The results of our systematic investigation with minimal hypotheses in \cref{res:generated_topology}, however, are more comprehensive than those in \cite{conradie:2005}.
\end{example}

\section{Hausdorff \uoLtops: going up and going down}\label{sec:uoLtops_going_up_and_going_down}

\noindent In this section, we investigate how, via a going-up-going-down construction, the existence of  a Hausdorff \oLtop\ on an order dense ideal of a vector lattice $E$ implies that every regular vector sublattice of $E$ admits a (necessarily unique) Hausdorff \uoLtop.

We start by going up.

 \begin{proposition}\label{res:preparation_uoLt}
 	Let $E$ be a vector lattice, and let $F$ be an ideal of $E$. Suppose that $F$ admits an \oLtop\ $\tau_F$. Choose a non-empty subset $S$ of $F$. Then $\unb_S\!\tau_F$ is a  \uoLtop\ on $E$. It is a \uppars{necessarily unique} Hausdorff \uoLtop\ on $E$ if and only if $\tau_F$ is a Hausdorff topology on $F$ and the ideal $I_S$ that is generated by $S$ is order dense in $E$.
 \end{proposition}

\begin{proof}
	We know from \cref{res:generated_topology} that $\unb_S\!\tau_F$ is a locally solid additive topology on $E$. In order to see that it is a linear topology on $E$, we verify the condition in \cref{eq:absorbing_condition}. Take $x$ in $E$ and $s$ in $S$. Then $\abs{\varepsilon x}\wedge \abs{s}\oconv 0$ in $E$ as $\varepsilon\to 0$ in $\RR$. Since $F$, being in ideal of $E$, is a regular vector sublattice of $E$, \cref{res:local-global_for_o-convergence_and_uo-convergence} shows that $\abs{\varepsilon x}\wedge\abs{s}\oconv 0$ in $F$. Since $\tau_F$ is an \oLtop\ on $F$, this implies that $\abs{\varepsilon x}\wedge\abs{s}\conv{\tau_F} 0$ in $F$ as $\varepsilon\to 0$ in $\RR$, as required.
	
	To conclude the proof, suppose that $\net$ is a net in $E$ such that $x_\alpha\uoconv 0$ in $E$. Take $s\in S$. Then $\abs{x_\alpha}\wedge \abs{s}\oconv 0$ in $E$. Again, since $F$ is a regular vector sublattice of $E$, \cref{res:local-global_for_o-convergence_and_uo-convergence} shows that $\abs{x_\alpha}\wedge \abs{s}\oconv 0$ in $F$. Since $\tau_F$ is an \oLtop\ on $F$, this implies that $\abs{x_\alpha}\wedge \abs{s}\conv{\tau_F} 0$ in $F$. It now follows from \cref{res:generated_topology} that $x_\alpha\conv{\unb_S\!\tau_F}0$ in $E$, as required.
	
	The uniqueness statement is clear from \cref{res:conradie_taylor}.
\end{proof}

The combination of \cref{res:generated_topology} and \cref{res:preparation_uoLt} immediately yields the following.

\begin{theorem}\label{res:Hausdorff_oLt_on_order_dense_ideal_generates_global_Hausdorff_uoLt}
Let $E$ be a vector lattice. Suppose that $E$ has an order dense ideal $F$ that admits a Hausdorff \oLtop\ $\tau_F$. Then $E$ admits a \uppars{necessarily unique} Hausdorff \uoLtop\ $\uoLt_E$. This topology $\uoLt_E$ is equal to $\unb_S\!\tau_F$ for every subset $S$ of $F$ such that the ideal $I_S\subseteq F$ that is generated by $S$ is order dense in $E$.

For a net $\net$ in $E$, the following are equivalent:
\begin{enumerate}
	\item $x_\alpha\conv{\uoLt_E} 0$ in $E$;\label{part:convergence_in_uoLtop}
	\item $\abs{x_\alpha}\wedge \abs{s}\conv{\tau_F}0$ in $F$ for all $s\in S$; \label{part:testing_against_subset}
	\item $\abs{x_\alpha}\wedge \abs{y}\conv{\tau_F}0$ in $F$ for all $y\in F$. \label{part:convergence_criteria_order_dense_ideal}
\end{enumerate}
\end{theorem}

\begin{remark}\label{rem:essential_difference}
For the case in \cref{res:Hausdorff_oLt_on_order_dense_ideal_generates_global_Hausdorff_uoLt} where $S=F$ and $\tau_F$ is the restriction of a Hausdorff \oLtop\ on $E$, it was already established in \cite[Theorem~9.6]{taylor:2019} that $\unb_F\!\tau_F$ is a Hausdorff \uoLtop\ on $E$. It is, therefore, of some importance to point out that not every Hausdorff \oLtop\ on an order dense ideal is the restriction of a Hausdorff \oLtop\ on the enveloping vector lattice. By way of example, consider the order dense ideal $c_0$ of $\ell^\infty$. Since the supremum norm on $c_0$ is order continuous, the usual norm topology $\tau_{c_0}$ on $c_0$ is a Hausdorff \oLtop. However, there does not even exist a possibly non-Hausdorff \oLtop\ $\tau_{\ell^\infty}$ on $\ell^\infty$ that extends $\tau_{c_0}$. In order to see this, consider the sequence of standard unit vectors $\seq{e}{n}$ in $\ell^\infty$. We have $e_n\oconv 0$ in $\ell^\infty$ , which would imply that $e_n\conv{\tau_{\ell^\infty}}0$ in $\ell^\infty$. Since $\tau_{\ell^\infty}$ extends $\tau_{c_0}$, we would have that $e_n\to 0$ in norm. This contradiction shows that such an extension does not exist.

Although the terminology is not used as such, the fact that $\unb_F\!\tau_F$ is a Hausdorff \uoLtop\ on $E$ is implicit in the construction in \cite[p.~290]{conradie:2005}.
\end{remark}

\begin{remark}\label{rem:small_subset_criterion}
    We are not aware of a reference where it is noted, as in part~\ref{part:testing_against_subset}, that convergence of a net in the Hausdorff \uoLtop\ on $E$ can be established by using a (presumably small and manageable) subset $S$ of $F$ instead of the full ideal $F$. This non-trivial fact, which relies on the uniqueness of a Hausdorff \uoLtop, appears to be of some practical value.
\end{remark}

In view of the uniqueness of a Hausdorff uo-Lebesgue topology (see \cref{res:conradie_taylor}), the following is now clear from \cref{res:Hausdorff_oLt_on_order_dense_ideal_generates_global_Hausdorff_uoLt}.

\begin{corollary}\label{res:same_generated_topology}
	Let $E$ be a vector lattice, and suppose that $E$ has order dense ideals $F_1$ and $F_2$, each of which admits a Hausdorff \oLtop. 	For $i=1,2$, choose a Hausdorff \oLtop\ $\tau_{F_i}$ on $F_i$, and choose a non-empty subset $S_i$ of $F_i$ such that the ideal $I_{S_i}\subseteq F_i$ that is generated by $S_i$ in $E$ is order dense in $E$. Then $\unb_{S_1}\!\tau_{F_1}$ and $\unb_{S_2}\!\tau_{F_2}$ are both equal to the \uppars{necessarily unique} \uoLtop\ $\uoLt_E$ on $E$.
\end{corollary}

\begin{remark}
	The case in \cref{res:same_generated_topology} where $S_1=F_1$ and $S_2=F_2$ is \cite[Proposition~3.2]{conradie:2005}.
	
The case where, for $i=1,2$, $S_i=F_i$ and $\tau_{F_i}$ is the restriction to $F_i$ of a Hausdorff \oLtop\ $\tau_i$ on $E$, is a part of \cite[Theorem~9.6]{taylor:2019}. Note, however, that our underlying proof in \cref{res:preparation_uoLt} that $\unb_S\!\tau_F$ is a \uoLtop\ is direct, whereas in the proof of \cite[Theorem~9.6]{taylor:2019} the identification of a Hausdorff \uoLtop\ as a minimal Hausdorff locally solid topology as in \cref{res:conradie_taylor} is used.
\end{remark}

Complementing the preceding going-up results, we cite the following going-down result; see \cite[Proposition 5.12]{taylor:2019}.

\begin{proposition}[Taylor]\label{res:going_down_one}
	Suppose that the vector lattice $E$ admits a \uppars{necessarily unique} Hausdorff \uoLtop\ $\uoLt_E$. Take a vector sublattice $F$ of $E$. Then $F$ is a regular vector sublattice of $E$ if and only if the restriction $\uoLt_E|_F$ of $\uoLt_E$ to $F$ is a \uppars{necessarily unique} Hausdorff \uoLtop\ on $F$.
\end{proposition}

A variation on this theme, with a wider range of topologies to use for testing the regularity of a vector sublattice, is the following.

\begin{proposition}\label{res:going_down_two}
		Suppose that the vector lattice $E$ admits a Hausdorff \oLtop\ $\tau_E$. Take a vector sublattice $F$ of $E$. Then $F$ is a regular vector sublattice of $E$ if and only if the restriction $\tau_E|_F$ of $\tau_E$ to $F$ is a Hausdorff \oLtop\ on $F$.
\end{proposition}

\begin{proof}
Once one recalls that, by definition, order convergence of a net to 0 in the regular vector sublattice $F$ of $E$ implies order convergence of the net to 0 in $E$, the proof is a straightforward minor adaptation of that of \cite[Proposition 5.12]{taylor:2019}.
\end{proof}

We now have the following overview theorem concerning Hausdorff \oLtops\ and Hausdorff \uoLtops\ on a vector lattice and on its order dense ideals. It is easily established by recalling that a \uoLtop\ is an \oLtop, that an ideal is a regular vector sublattice, and by using \cref{res:Hausdorff_oLt_on_order_dense_ideal_generates_global_Hausdorff_uoLt}, \cref{res:going_down_one}, and \cref{res:going_down_two}.

\begin{theorem}\label{res:overview}
	Let $E$ be a vector lattice, and let $F$ be an order dense ideal of $E$.
	\begin{enumerate}
		\item Suppose that $E$ admits a Hausdorff \oLtop\ $\tau_E$. Then the restricted topology $\tau_E|_F$ is a Hausdorff \oLtop\ on $E$.
		\item Suppose that $E$ admits a \uppars{necessarily unique} Hausdorff \uoLtop\ $\uoLt_E$. Then the restricted topology $\uoLt_E|_F$ is a \uppars{necessarily unique} Hausdorff \uoLtop\ on $F$.
		\item The following are equivalent:
		\begin{enumerate}
			\item $F$ admits a Hausdorff \oLtop;
			\item $F$ admits a Hausdorff \uoLtop;
			\item $E$ admits a Hausdorff \oLtop;
			\item $E$ admits a Hausdorff \uoLtop.
		\end{enumerate}	
			In that case, the unique \uoLtop\ $\uoLt_E$ on $E$ equals $\unb_S\!\tau_F$ for every Hausdorff \oLtop\ on $F$ and every subset $S$ of $F$ such that the ideal $I_S\subseteq F$ is order dense in $E$, and the following are equivalent:
				\begin{enumerate_roman}
					\item $x_\alpha\conv{\uoLt_E} 0$ in $E$;
					\item $\abs{x_\alpha}\wedge \abs{s}\conv{\tau_F}0$ in $F$ for all $s\in S$;
					\item $\abs{x_\alpha}\wedge \abs{y}\conv{\tau_F}0$ in $F$ for all $y\in F$.
				\end{enumerate_roman}		
	\end{enumerate}
\end{theorem}

We conclude this section with a short discussion of Banach lattices with order continuous norms. Evidently, the norm topologies on such Banach lattices are Hausdorff \oLtops. As already noted in \cite[p.~993]{taylor:2019}, \cref{res:Hausdorff_oLt_on_order_dense_ideal_generates_global_Hausdorff_uoLt} allows one to identify the so-called un-topologies (see \cite[Section~7]{deng_o_brien_troitsky:2017} and \cite[p.~746]{kandic_li_troitsky:2018}) on such lattices as the Hausdorff \uoLtops\ that these spaces apparently admit. Consequently, we have the following result. The case where $S=E$ can be found in \cite[p.~993]{taylor:2019}.

\begin{proposition}\label{res:banach_lattice_with_order_continuous_norm_one}
	Let $E$ be a Banach lattice with an order continuous norm and norm topology $\tau_E$. Then $E$ admits a \uppars{necessarily unique} \uoLtop.
	
	Choose a subset $S$ of $E$ such that the ideal $I_S$ that is generated by $S$ in $E$ is order dense in $E$. Then:
	\begin{enumerate}
		\item $\unb_S\!\tau_E$ is the uo-Lebesgue topology $\uoLt_E$ of $E$;
		\item when $\net$ is a net in $E$, then $x_\alpha\conv{\uoLt_E} 0$ in $E$ if and only if $\norm{\abs{x_\alpha}\wedge \abs{s}}\conv{}0$ for all $s\in S$; equivalently, if and only if $\norm{\abs{x_\alpha}\wedge \abs{y}}\conv{}0$ for all $y\in E$.\label{part:convergence_criteria_order_dense_ideal_banac_lattice_with_oc_norm}
	\end{enumerate}
\end{proposition}

There is an alternative reason why Banach lattices with an order continuous norms admit Hausdorff \uoLtops, and this  results in an alternative description of these topologies; see \cref{res:banach_lattices_with_order_continuous_norm_two}, below.

Finally, suppose that $E$ is a vector lattice that has order dense ideals $F_1$ and $F_2$ that are Banach lattices with order continuous norm topologies $\tau_{F_1}$ and  $\tau_{F_2}$, respectively. Then it is immediate from \cref{res:same_generated_topology} that $E$ admits a Hausdorff \uoLtop\ $\uoLt_E$, and that $\unb_{F_1}\!\tau_{F_1}$ and $\unb_{F_2}\!\tau_{F_2}$ are both equal to $\uoLt_E$. As discussed in \cref{exam:kandic_li_troitsky:2018}, this can, using the terminology in \cite{kandic_li_troitsky:2018}, be rephrased as stating that $F_1$ and $F_2$ induce the same un-topology on $E$. We have thus retrieved \cite[Theorem~2.6]{kandic_li_troitsky:2018}.

\section{\uoLtops\ generated by absolute weak topologies on order dense ideals}\label{sec:separating_order_continuous_dual}

\noindent In this section, we shall be concerned with vector lattices having order dense ideals with separating order continuous duals as a source for Hausdorff \uoLtops\ on the vector lattices themselves.

We start by recapitulating some facts from \cite[p.~63--64]{aliprantis_burkinshaw_LOCALLY_SOLID_RIESZ_SPACES_WITH_APPLICATIONS_TO_ECONOMICS_SECOND_EDITION:2003}. Let $E$ be a vector lattice, and let $A$ be a non-empty subset of the order dual $\odual{E}$ of $E$. For $\varphi\in A$, define the lattice semi-norm $\rho_\varphi: E\to[0,\infty)$ by setting $\rho_\varphi(x)\coloneqq \abs{\varphi}(\abs{x})$ for $x\in E$.
Then the locally convex-solid linear topology on $E$ that is generated by the family $\{\rho_\varphi: \varphi\in A\}$ is called the \emph{absolute weak topology generated by $A$ on $E$}; it is denoted by $\abs{\sigma}(E,A)$. With $I_A$ denoting the ideal generated by $A$ in $\odual{E}$, we have $\abs{\sigma}(E,A)=\abs{\sigma}(E,I_A)$. Using \cref{res:polars}, one easily concludes that $\abs{\sigma}(E,A)$ is Hausdorff if and only if $I_A$ separates the points of $E$. Although we shall not use it, let us still remark that it is not  difficult to see that a net $\net$ in $E$ is $\abs{\sigma}(E,A)$-convergent to zero if and only if $\varphi(x_\alpha)\conv{} 0$ uniformly for $\varphi$ in each fixed order interval of $I_A$. Thus absolute weak topologies are more natural than is perhaps apparent from their definition.

The following is now clear.

\begin{lemma}\label{res:absolute_weak_topology_has_the_o-lebesgue_property}
	Let $E$ be a vector lattice, and let $A$ be a non-empty subset of $\ocdual{E}$. Let $I_A$ denote the ideal that is generated by $A$ in $\ocdual{E}$. Then $\abs{\sigma}(E,A)=\abs{\sigma}(E,I_A)$ is an \oLtop\ on $E$ that is even locally convex-solid. It is a Hausdorff topology if and only if $I_A$ separates the points of $E$. When $\net$ is a net in $E$, then $x_\alpha\conv{\abs{\sigma}(E,A)}0$ in $E$ if and only if $\abs{\varphi}\left(\abs{x_\alpha}\right)\to 0$ for all $\varphi\in A$; equivalently, if and only if $\abs{\varphi}\left(\abs{x_\alpha}\right)\to 0$ for all $\varphi\in I_A$.
\end{lemma}

Now that \cref{res:absolute_weak_topology_has_the_o-lebesgue_property} provides a whole class of vector lattices admitting Hausdorff \oLtops, we can use these as input for \cref{res:Hausdorff_oLt_on_order_dense_ideal_generates_global_Hausdorff_uoLt}. Taking the convergence statements in \cref{res:absolute_weak_topology_has_the_o-lebesgue_property} into account, we arrive at the following.

\begin{theorem}\label{res:uo-lebesgue_topology_originating_from_separating_order_continuous_dual_of_order_dense_ideal}
	Let $E$ be a vector lattice. Suppose that $E$ has an order dense ideal $F$ such that $\ocdual{F}$ separates the points of $F$. Then $E$ admits a \uppars{necessarily unique} Hausdorff \uoLtop\ $\uoLt_E$.
	
	Choose a subset $A$ of $\ocdual{F}$ such that the ideal $I_A$ that is generated by $A$ in $\ocdual{F}$ separates the points of $F$, and choose a subset $S$ of $F$ such that the ideal $I_S\subseteq F$ that is generated by $S$ is order dense in $E$. Then:
	\begin{enumerate}
		\item $\unb_S\abs{\sigma}(F,A)$ and $\unb_F\abs{\sigma}(F,I_A)$ are both equal to $\uoLt_E$;
		\item for a net $\net$ in $E$, $x_\alpha\conv{\uoLt_E}0$ in $E$ if and only if $\abs{\varphi}\left(\abs{x_\alpha}\wedge\abs{s}\right)\to 0$ for all $\varphi\in A$ and $s\in S$; equivalently, if and only if $\abs{\varphi}\left(\abs{x_\alpha}\wedge\abs{y}\right)\to 0$ for all $\varphi\in \ocdual{F}$ and $y\in F$.
	\end{enumerate}
\end{theorem}

For the sake of completeness, we recall that a regular vector sublattice of a vector lattice $E$ as in the theorem also has a (necessarily unique) Hausdorff \uoLtop, and that this topology is the restriction of $\uoLt_E$ to the vector sublattice.

\begin{remark}
As noted in \cref{rem:explicit_neighbourhood_base}, one can give an explicit neighbourhood base at zero for the topology $\uoLt_E$ in \cref{res:uo-lebesgue_topology_originating_from_separating_order_continuous_dual_of_order_dense_ideal}.
\end{remark}

For Banach lattices with order continuous norms, the order/norm dual consists of order continuous linear functionals only. Hence we have the following consequence of \cref{res:uo-lebesgue_topology_originating_from_separating_order_continuous_dual_of_order_dense_ideal}.

\begin{corollary}\label{res:banach_lattices_with_order_continuous_norm_two}
A Banach lattice $E$ with an order continuous norm admits a \uppars{necessarily unique} Hausdorff \uoLtop\ $\uoLt_E$, namely $\unb_E\abs{\sigma}(E,E^\ast)$.	
\end{corollary}	

\begin{remark}
For a Banach lattice $E$ with an order continuous norm, \cref{res:banach_lattice_with_order_continuous_norm_one} also shows that the (necessarily unique) Hausdorff \uoLtop\ $\uoLt_E$ on $E$ is the unbounded norm topology. The unbounded norm topology and the unbounded absolute weak topology therefore coincide for Banach lattices with order continuous norms.
\end{remark}

The following gives a necessary condition for convergence in a Hausdorff \uoLtop. It is essential in the proof of \cref{res:tau_m_to_sub_uo}, below.

\begin{proposition}\label{res:essential}
 Let $E$ be a vector lattice that admits a \uppars{necessarily unique} Hausdorff \uoLtop\ $\uoLt_E$, and let  $\net$ be a net in $E$ such that $x_\alpha\conv{\uoLt_E}0$ in $E$. Take an ideal $F$ of $E$ such that $\ocdual{F}$ separates the points of $F$.  Then $\abs{\varphi}(\abs{x_\alpha}\wedge \abs{y})\conv{}0$ for all $\varphi\in \ocdual{F}$ and $y\in F$.
\end{proposition}

\begin{proof}
	Take $\varphi\in\ocdual{F}$ and $y\in F$. Since $\uoLt_E$ is a locally solid topology, we have $\abs{x}_\alpha\wedge\abs{y}\conv{\uoLt_E}0$ in $E$. It follows from \cref{res:going_down_one} that $F$ has a (necessarily unique) Hausdorff \uoLtop\ $\uoLt_F$ and that $\abs{x}_\alpha\wedge\abs{y}\conv{\uoLt_F}0$. Now we apply \cref{res:uo-lebesgue_topology_originating_from_separating_order_continuous_dual_of_order_dense_ideal} with $E=F$ to see that $\abs{\varphi}((\abs{x_\alpha}\wedge \abs{y})\wedge\abs{y})\conv{}0$.
\end{proof}	

We shall now consider the order dual $\odual{E}$ of a vector lattice $E$. For $x\in E$, we set
\[
\varphi_x(\varphi)\coloneqq \varphi(x)
\]
for $\varphi\in \odual{E}$. Then $\varphi_x\in\ocdual{\left(\odual{E}\right)}$, and the map $\varphi:E\to \odual{E}$ is a lattice homomorphism; see \cite[p.~43]{aliprantis_burkinshaw_LOCALLY_SOLID_RIESZ_SPACES_WITH_APPLICATIONS_TO_ECONOMICS_SECOND_EDITION:2003}. Since $\varphi(E)$ already separates the points of $\odual{E}$, we see that $\ocdual{\left(\odual{E}\right)}$ separates the points of $\odual{E}$.

We can now apply \cref{res:uo-lebesgue_topology_originating_from_separating_order_continuous_dual_of_order_dense_ideal} twice. In both cases, we replace $E$ with $\odual{E}$, and we choose $\odual{E}$ for both $F$ and $S$. In the first application, we choose $\ocdual{\left(\odual{E}\right)}$ for $A$; in the second, we choose $\varphi(E)$. The result is as follows.

\begin{corollary}\label{res:corollary_to_uniqueness_of_tau_m_for_dual}
Let $E$ be a vector lattice. Then the order dual  $\odual{E}$ of $E$ admits a \uppars{necessarily unique} Hausdorff \uoLtop\ $\uoLt_{\odual{E}}$.

Moreover:
\begin{enumerate}
	\item $\unb_{\odual{E}}\abs{\sigma}( \odual{E},  \ocdual{\left(\odual{E}\right)})$ and $\unb_{\odual{E}}\abs{\sigma}( \odual{E},  E)$ are both equal to $\uoLt_{\odual{E}}$;
	\item when $(\varphi_\alpha)_{\alpha\in\is}$ is a net in $\odual{E}$, then we have that  $\varphi_\alpha\conv{\uoLt_{\odual{E}}}0$ in $E$ if and only if $\abs{\xi}\left(\abs{\varphi_\alpha}\wedge\abs{\varphi}\right)\to 0$ for all  $\xi\in\ocdual{\left(\odual{E}\right)}$ and $\varphi\in \odual{E}$; equivalently, if and only if $(\abs{\varphi_\alpha}\wedge\abs{\varphi})(\abs{x})\to 0$ for all  $x\in E$ and $\varphi\in \odual{E}$.
\end{enumerate}
\end{corollary}

\begin{remark}\quad
	\begin{enumerate}
	\item 	As in the case of \cref{res:uo-lebesgue_topology_originating_from_separating_order_continuous_dual_of_order_dense_ideal},  \cref{rem:explicit_neighbourhood_base} shows how to give an explicit neighbourhood base at zero for the topology $\uoLt_{\odual{E}}$ in \cref{res:corollary_to_uniqueness_of_tau_m_for_dual}.
	\item By \cref{res:going_down_one}, every regular sublattice of the order dual of a vector lattice also admits a (necessarily unique) Hausdorff Lebesgue topology that can be described in two ways. For an ideal, one of these descriptions is already in \cite[Example~5.8]{taylor:2019}.
\item \cref{res:corollary_to_uniqueness_of_tau_m_for_dual} shows that, in particular, the norm/order dual $E^\ast$ of a Banach lattice admits a (necessarily unique) Hausdorff \uoLtop\ $\uoLt_{E^\ast}$, namely the so-called unbounded absolute weak $^\ast$-topology $\unb_{E^\ast}\abs{\sigma}(E^\ast,E)$.  This was already observed in \cite[Lemma~6.6]{taylor:2019}.
	\end{enumerate}
\end{remark}

\section{Regular vector sublattices of $\Ell_0(X,\Sigma,\mu)$ for semi-finite measures}\label{sec:vector_lattices_of_equivalence_classes_of_measurable_functions}

\noindent Let $(X,\Sigma,\mu)$ be a measure space, and write $\Ell_0(X,\Sigma,\mu)$ for the vector lattice of all real-valued $\Sigma$-measurable functions on $X$, with identification of two functions when they agree $\mu$-almost everywhere. In this section we show that, for semi-finite $\mu$, every regular sublattice of $\Ell_0(X,\Sigma,\mu)$ admits a (necessarily unique) Hausdorff \uoLtop, and that a net converges in this topology if and only if it converges in measure on subsets of finite measure; see \cref{res:going_downlattice_of_L_0}, below.

For some regular sublattices of $\Ell_0(X,\Sigma,\mu)$, it is quite obvious that they admit a Hausdorff \uoLtop. Recall that the spaces $\Ell_p(X,\Sigma,\mu)$ for $p$ such that $1\leq p<\infty$ have order continuous norms for all measures $\mu$; see  \cite[Theorem~13.7]{aliprantis_border_INFINITE_DIMENSIONAL_ANALYSIS_THIRD_EDITION:2006}, for example. Hence their norm topologies are Hausdorff \oLtops, and then their un-topologies are the Hausdorff \uoLtops\ on these spaces. Alternatively, one can observe that their order continuous duals separate their points, and then also identify the Hausdorff \uoLtops\ on these spaces as the unbounded absolute weak topologies. In a similar vein, when $\mu$ is $\sigma$-finite, every ideal of $L_0(X,\Sigma,\mu)$ that can be supplied with a lattice norm has a separating order continuous dual. This result of Lozanovsky's (see \cite[Theorem~5.25]{abramovich_aliprantis_INVITATION_TO_OPERATOR_THEORY:2002}, for example) then implies that such a normed function space admits a Hausdorff \uoLtop.

How about the spaces $\Ell_p(X,\Sigma,\mu)$ for $0\leq p<1$? There is no norm to work with, and it may well be the case that their order continuous duals are even trivial. Indeed, when $\mu$ is atomless, then, according to a results of Day's, the order continuous dual of $\Ell_p(X,\Sigma,\mu)$ is trivial for $0<p<1$; see \cite[Theorem~13.31]{aliprantis_border_INFINITE_DIMENSIONAL_ANALYSIS_THIRD_EDITION:2006}, for example. According to \cite[Exercise~25.2]{zaanen_INTRODUCTION_TO_OPERATOR_THEORY_IN_RIESZ_SPACES:1997}, the order continuous dual of $\Ell_0(X,\Sigma,\mu)$ is trivial for every $\sigma$-finite measure with the property that, for any measurable subset $A$ such that $0<\mu(A)<\infty$ and for any $\alpha$ such that $0<\alpha<\mu(A)$, there exists a measurable subset $A^\prime$ of $A$ such that $\mu(A^\prime)=\alpha$. Taking \cite[Exercise~10.12 on p.~67]{zaanen_INTEGRATION:1967} into account, we see that, in particular, the order continuous dual of $\Ell_0(X,\Sigma,\mu)$ is trivial for all atomless $\sigma$-finite measures.

In spite of the failure of the two obvious approaches, it is still possible to show that all spaces $\Ell_p(X,\Sigma,\mu)$ for $0\leq p<1$ admit Hausdorff \uoLtops, provided that the measure is semi-finite. For such $\mu$, this is even true for all regular vector sublattices of $\Ell_0(X,\Sigma,\mu)$. This can be seen via the going-up-going-down approach from \cref{sec:uoLtops_going_up_and_going_down}, and we shall now elaborate on this. We start with a few preliminary remarks.

Recall that a measure space $(X, \Sigma, \mu)$ is said to be \emph{semi-finite} if, for any $A\in \Sigma$ with $\mu(A)=\infty$, there exists a measurable subset $A^\prime$ of $A$ such that $0<\mu(A^\prime)<\infty$. Every $\sigma$-finite measure is semi-finite. For an arbitrary measure $\mu$ and an arbitrary $p$ such that $1\leq p<\infty$, it is easy to see that the ideal $\Ell_p(X,\Sigma,\mu)$ of $\Ell_0(X,\Sigma,\mu)$ is order dense in $\Ell_0(X,\Sigma,\mu)$ if and only if $\mu$ is semi-finite.  In that case, the ideal that is generated in $\Ell_0(X,\Sigma,\mu)$ by the subset $S\coloneqq \{1_A: A\in\Sigma\text{ has finite measure}\}$ of $\Ell_p(X,\Sigma,\mu)$ is obviously also order dense in $\Ell_0(X,\Sigma,\mu)$.

Let $(X,\Sigma,\mu)$ be a measure space.  Take $f\in\Ell_0(X,\Sigma,\mu)$. Then a net $(f_\alpha)_{\alpha\in\is}$ in $\Ell_0(X,\Sigma,\mu)$ \emph{converges to f in measure on subsets of finite measure} when, for all $A\in\Sigma$ such that $\mu(A)<\infty$ and for all $\varepsilon>0$,  $\mu(\{x\in A: \abs{f_\alpha(x)-f(x)}\geq\varepsilon\})\conv{}0$.  In that case, we write $f_\alpha\conv{{\mu^\ast}}f$, using as asterisk to distinguish this convergence from the perhaps more usual global convergence in measure.

The following is the core result of this section. We recall that, as already mentioned, the spaces $\Ell_p(X,\Sigma,\mu)$ have order continuous norms for all measures $\mu$ and for all $p$ such that $1\leq p<\infty$, so that their norm topologies are Hausdorff \oLtops.

\begin{theorem}\label{res:tau_m_is_convergence_in_measure}
	Let $E=\Ell_0(X,\Sigma,\mu)$, where $\mu$ is a semi-finite measure. Then $G$ admits a \uppars{necessarily unique} Hausdorff \uoLtop\ $\uoLt_E$.
	
	Take a net $(f_\alpha)_{\alpha\in\is}$ in $E$. Then the following are equivalent for every $p$ such that $1\leq p<\infty$:
	\begin{enumerate}
		\item $f_\alpha\conv{\uoLt_E} 0$; \label{part:Ell_0_1}
		\item
		\[
		\int_X\! \abs{f_\alpha}^p\wedge 1_A\,\text{d}\mu=\norm{\,\abs{f_\alpha}\wedge \abs{1_A}\,}_p^p\conv{}0
		\]
		for every measurable subset $A$ of $X$ with finite measure;\label{part:Ell_0_2}
		\item
		\[
		\int_X\! \abs{f_\alpha}^p\wedge \abs{f}^p\,\text{d}\mu=\norm{\,\abs{f_\alpha}\wedge \abs{f}\,}_p^p\conv{}0
		\]for every $f\in\Ell_p(X,\Sigma,\mu)$;\label{part:Ell_0_3}
		\item $f_\alpha\conv{\mu^\ast} f$.\label{part:Ell_0_4}
	\end{enumerate}
\end{theorem}

\begin{proof}
	
	We know from the semi-finiteness of $\mu$ that, for $p$ such that $1\leq p\leq \infty$, $L_p(X,\Sigma,\mu)$ is an order dense ideal of $\Ell_0(X,\Sigma,\mu)$.  Since $L_p(X,\Sigma,\mu)$ admits a Hausdorff \oLtop\ when $1\leq p<\infty$, \cref{res:Hausdorff_oLt_on_order_dense_ideal_generates_global_Hausdorff_uoLt} shows that $\Ell_0(X,\Sigma,\mu)$ admits a (necessarily unique) Hausdorff \uoLtop, and also that the statements in the parts~\ref{part:Ell_0_1},~\ref{part:Ell_0_2}, and~\ref{part:Ell_0_3} of the present theorem are equivalent for all such $p$.
	
	We show that part~\ref{part:Ell_0_3} implies part~\ref{part:Ell_0_4}.
	Take a measurable subset $A$ of $X$ with finite measure, and let $\varepsilon>0$. Since $\varepsilon 1_A\in\Ell_p(X,\Sigma,\mu)$, we have, by assumption,
	\[
	\int_X\! \abs{f_\alpha}^p\wedge(\varepsilon^p 1_A)\,\text{d}\mu\conv{}0.
	\]
	Because
	\[
	\int_X \!\abs{f_\alpha}^p\wedge(\varepsilon^p 1_A)\,\text{d}\mu\geq \int_{\{x\in A: \abs{f_\alpha(x)}\geq\varepsilon\}}\!\!\! \varepsilon^p\,\text{d}\mu=\varepsilon^p\mu\left(\{x\in A: \abs{f_\alpha(x)}\geq\varepsilon\}\right)
	\]
	we conclude that $\mu(\{x\in A: \abs{f_\alpha(x)}\geq\varepsilon\})\conv{}0$. Hence $f_\alpha\conv{\mu^\ast} 0$.
	
	We show that part~\ref{part:Ell_0_4} implies part~\ref{part:Ell_0_2}. Take a measurable subset $A$ of $X$ with finite measure, and take $\varepsilon>0$. Choose a $\delta>0$ such that $\delta^p\mu(A)<\varepsilon/2$. Then
	\begin{align*}
	\int_X\!  \abs{f_\alpha}^p\wedge1_A\,\text{d}\mu&= \int_{\{x\in A: \abs{f_\alpha(x)}^p\geq\delta^p\}}\!\!\!\abs{f_\alpha}^p\wedge 1_A \,\text{d}\mu+\int_ {\{x\in A: \abs{f_\alpha(x)}^p<\delta^p\}}\!\!\!\abs{f_\alpha}^p\wedge 1_A\,\text{d}\mu\\
	&\leq \int_{\{x\in A: \abs{f_\alpha(x)}^p\geq\delta^p\}}\!\!\!1 \,\text{d}\mu + \int_A\delta^p\,\text{d}\mu\\
	&\leq \mu\left(\{x\in A: \abs{f_\alpha(x)}\geq\delta\}\right)+ \varepsilon/2.
	\end{align*}
	By our  assumption, there exists an $\alpha_0\in\is$ such that  $\mu\left(\{x\in A: \abs{f_\alpha(x)}\geq\delta\}\right)<\varepsilon/2$ for all $\alpha\geq\alpha_0$. Then $\int_X  \abs{f_\alpha}^p\wedge1_A\,\text{d}\mu<\varepsilon$ for all $\alpha\geq\alpha_0$. Hence $\int_X \abs{f_\alpha}\wedge1_A\,\text{d}\mu\conv{}0$.
\end{proof}

\begin{remark}\quad
	\begin{enumerate}
		\item We are not aware of a proof of \cref{res:tau_m_is_convergence_in_measure} in the literature.
It is stated in \cite[p.~292]{conradie:2005} that the parts~\ref{part:Ell_0_1} and~\ref{part:Ell_0_4} are equivalent, but there only a reference is given to \cite[65K and~63L]{fremlin_TOPOLOGICAL_RIESZ_SPACES_AND_MEASURE_THEORY:1974}. Since  \cite[63L]{fremlin_TOPOLOGICAL_RIESZ_SPACES_AND_MEASURE_THEORY:1974} relies on the solution of the non-trivial exercise \cite[Exercise~63M(j)]{fremlin_TOPOLOGICAL_RIESZ_SPACES_AND_MEASURE_THEORY:1974} for which a solution is not provided, we thought it appropriate to give an independent proof in the present paper.
		\item The equivalence of the parts~\ref{part:Ell_0_3} and~\ref{part:Ell_0_4} for finite measures and sequences was also established by different methods in \cite[Example~23]{troitsky:2004}. Still earlier, this case was covered in \cite[Corollary~4.2]{deng_o_brien_troitsky:2017}, with a proof in the same spirit as our proof.
	\end{enumerate}
\end{remark}

As an immediate consequence of \cref{res:going_down_one} and \cref{res:tau_m_is_convergence_in_measure}, we obtain the following result via our going-up-going-down approach.

\begin{theorem}\label{res:going_downlattice_of_L_0}
	Let $(X,\Sigma,\mu)$ be a measure space, where $\mu$ is a semi-finite measure. Take a regular vector sublattice $E$ of $\Ell_0(X,\Sigma,\mu)$. Then $E$ admits a \uppars{necessarily unique} Hausdorff \uoLtop\ $\uoLt_E$. This
	topology $\uoLt_E$ on $E$ is the restriction of the Hausdorff \uoLtop\ on $\Ell_0(X,\Sigma,\mu)$. A net $(f_\alpha)_{\alpha\in\is}$ in $E$  converges to zero in $\uoLt_E$ if and only if it satisfies one of the three equivalent criteria in the parts~\ref{part:Ell_0_2},~\ref{part:Ell_0_3}, and~\ref{part:Ell_0_4} of \cref{res:tau_m_is_convergence_in_measure}. In particular, it is $\uoLt_E$-convergent to zero if and only if it converges to zero in measure on subsets of finite measure.
\end{theorem}

\begin{remark}\label{rem:various_descriptions}
	Let $p$ be such that $1\leq p<\infty$. For arbitrary measures, \cref{res:banach_lattice_with_order_continuous_norm_one,res:banach_lattices_with_order_continuous_norm_two} both give a description of the convergent nets in the Hausdorff \uoLtop\ on $\Ell_p(X,\Sigma,\mu)$. The former as the convergent nets in the un-topology, and the latter as the convergent nets in the unbounded absolute weak topology, respectively. When $\mu$ is semi-finite, \cref{res:going_downlattice_of_L_0} gives a third description as the convergence in measure on subsets of finite measure.
	
	Also for $p=\infty$, \cref{res:going_downlattice_of_L_0} shows that $\Ell_\infty(X,\Sigma,\mu)$ admits a (necessarily unique) Hausdorff \uoLtop\ whenever $\mu$ is semi-finite, and gives a description of its convergent nets. When $\mu$ is a localisable measure, two more descriptions are possible. We refer to \cite[211G]{fremlin_MEASURE_THEORY_VOLUME_2:2003} for the definition of localisable measures, and note that $\sigma$-finite measures are localisable, and that localisable measures are semi-finite. Indeed, for localisable measures, $\Ell_\infty(X,\Sigma,\mu)$ is the order dual of $L_1(X,\Sigma,\mu)$; see  \cite[243G(b)]{fremlin_MEASURE_THEORY_VOLUME_2:2003}.  Hence \cref{res:corollary_to_uniqueness_of_tau_m_for_dual} shows once more that $\Ell_\infty(X,\Sigma,\mu)$ admits a Hausdorff \uoLtop\ when $\mu$ is localisable, and gives a second and third description of its convergent nets.
\end{remark}

\begin{remark}
	Let $(X,\Sigma,\mu)$ be a measure space, where $\mu$ is a semi-finite measure.
	
	Let $p$ be such that $0< p<\infty$. On combining \cref{res:going_downlattice_of_L_0} and \cref{rem:minimal_and_smallest}, we see that the topology of convergence in measure on subsets of finite measure is the \emph{smallest} Hausdorff locally solid linear topology on $\Ell_p(X,\Sigma,\mu)$.\footnote{For this conclusion, we should note here that the usual metric topology on $\Ell_p(X,\Sigma,\mu)$ is a complete \oLtop\ for every measure $\mu$ and for every $p$ such that $0< p<\infty$. This is commonly known when $1\leq p<\infty$. When $0<p<1$, then the completeness is asserted in \cite[1.47]{rudin_FUNCTIONAL_ANALYSIS_SECOND_EDITION:1991}. The fact that the metric topology is an \oLtop\ for such $p$ follows from what is stated on \cite[p.~211]{aliprantis_burkinshaw_LOCALLY_SOLID_RIESZ_SPACES_WITH_APPLICATIONS_TO_ECONOMICS_SECOND_EDITION:2003} in the context of $\sigma$-finite measures. This implies the result for general measures. Indeed, suppose that $\netgen{f_\alpha}{\alpha\in \is}$ is a net in $\Ell_p(X,\Sigma,\mu)$ such that $f_\alpha\downarrow 0$. Passing to a tail, we may suppose that the net is bounded above by an $f_{\alpha_0}\in\Ell_p(X,\Sigma,\mu)$. The support of this $f_{\alpha_0}$ is $\sigma$-finite. Using the fact that the elements of $\Ell_p(X,\Sigma,\mu)$ that vanish off this support form an ideal of $\Ell_p(X,\Sigma,\mu)$, it is then easily seen from the $\sigma$-finite case that the chosen tail of the net converges to zero in the metric topology of $\Ell_p(X,\Sigma,\mu)$. }    For $\sigma$-finite measures, this can already be found in \cite[Theorem~7.74]{aliprantis_burkinshaw_LOCALLY_SOLID_RIESZ_SPACES_WITH_APPLICATIONS_TO_ECONOMICS_SECOND_EDITION:2003}, where it is also established that the usual metric topology is then the largest Hausdorff locally solid linear topology.
	
	For $p=\infty$, the combination of \cref{res:going_downlattice_of_L_0} and \cref{res:conradie_taylor} shows that the topology of convergence in measure on subsets of finite measure is the unique \emph{minimal} Hausdorff locally solid linear topology on $\Ell_\infty(X,\Sigma,\mu)$. It seems worthwhile to note that, when $\mu$ is, in fact, $\sigma$-finite, and also non-atomic, \cite[Theorem~7.75]{aliprantis_burkinshaw_LOCALLY_SOLID_RIESZ_SPACES_WITH_APPLICATIONS_TO_ECONOMICS_SECOND_EDITION:2003} shows that there is now no \emph{smallest} Hausdorff locally solid linear topology on $\Ell_\infty(X,\Sigma,\mu)$.
\end{remark}

\begin{remark}
	Let $\seq{x_n}{n}$ be a sequence in $\Ell_0(X,\Sigma,\mu)$, where $\mu$ is a semi-finite measure. Suppose that $f_n\conv{}0$ $\mu$-almost everywhere. Then $f_n\conv{\mu^\ast}0$. This is immediate from Egoroff's theorem (see \cite[Theorem~2.33]{folland_REAL_ANALYSIS_SECOND_EDITION:1999}, for example), but it can also be obtained (with a long detour) in the context of uo-convergence and \uoLtops. Indeed, by \cite[Proposition~3.1]{gao_troitsky_xanthos:2017}, almost everywhere convergence of a sequence in $\Ell_0(X,\Sigma,\mu)$ is, for arbitrary measures, equivalent to uo-convergence in $\Ell_0(X,\Sigma,\mu)$. Since, by definition, uo-convergence implies convergence in a \uoLtop\ (when this exists), an appeal to \cref{res:tau_m_is_convergence_in_measure} also yields the desired result.
\end{remark}

\section{uo-convergent sequences within $\uoLt_E$-convergent nets}\label{sec:uo-convergent_subsequences_of_uoLt-convergent_nets}

\noindent Let $E$ be a vector lattice that admits a (necessarily unique) Hausdorff \uoLtop\ $\uoLt_E$. When $\net$ is a net in $E$ such that $x_\alpha\uoconv 0$, then, by definition, $x_\alpha\conv{\uoLt_E} 0$. The present section is concerned with results that go in the opposite direction. The main result is \cref{res:tau_m_to_sub_uo}, below, which lies at the basis of topological considerations in \cref{sec:topological aspects of uo-convergence}, but we start with a few more elementary results.

For an atomic vector lattice $E$, the situation is as easy as can be.  Recall that, by \cite[Theorem~1.78]{aliprantis_burkinshaw_LOCALLY_SOLID_RIESZ_SPACES_WITH_APPLICATIONS_TO_ECONOMICS_SECOND_EDITION:2003}, the atomic vector lattices are precisely the order dense vector sublattices of $\RR^X$ for some set $X$. Combining \cite[Lemma~3.1]{dabboorasad_emelyanov_marabeh:2020} and \cite[Lemma~7.4]{taylor:2019}, we have the following.

\begin{proposition}[Taylor]\label{res:atomic}
	Let $E$ be an atomic vector lattice. Then $E$ admits a \uppars{necessarily unique} Hausdorff \uoLtop\ $\uoLt_E$, and this topology is locally convex-solid. For a net in $E$, uo-convergence and $\uoLt_E$-convergence coincide, so that uo-convergence is topological. When $E$ is an order dense vector sublattice of $\RR^X$ for some set $X$, then a net in $E$ is uo- and $\uoLt_E$-convergent if and only if it is pointwise convergent.
\end{proposition}

For monotone nets, uo-convergence and $\uoLt_E$-convergence both coincide with order convergence, according to the following lemma that is a consequence of  \cite[Theorem~2.21]{aliprantis_burkinshaw_LOCALLY_SOLID_RIESZ_SPACES_WITH_APPLICATIONS_TO_ECONOMICS_SECOND_EDITION:2003}.

\begin{lemma}\label{res:monotone_tau_and_uo}
	Let $E$ be a vector lattice, and suppose that $\tau$ is a Hausdorff locally solid linear topology on $E$.  Let $\net$ be a monotone net in $E$ and let $x\in E$. When $x_\alpha\tauconv x$ in $E$, then $x_\alpha\oconv x$ in $E$, which is equivalent to $x_\alpha\uoconv x$. When $\uoLt_E$ is a \uppars{necessarily unique} Hausdorff \uoLtop\ on $E$, then $x_\alpha\conv{\uoLt_E} x$ in $E$ if and only if $x_\alpha\oconv x$ in $E$ if and only if $x_\alpha\uoconv x$.
\end{lemma}

For non-monotone nets in general vector lattices, it is not generally true that $\uoLt_E$-convergence implies uo-convergence. This can already fail for sequences in Banach lattices with order continuous norms. As an example, consider $E=\Ell_1([0,1])$. For $n=1,2,\dotsc$ and $k=1,2,\dotsc,n$, let $f_{nk}$ be the characteristic function of $[\frac{k-1}{n},\frac{k}{n}]$, and consider the sequence $f_{11}, f_{21}, f_{22}, f_{31}, f_{32}, f_{33}, f_{41},\dotsc$. It converges to zero in measure, so \cref{res:tau_m_is_convergence_in_measure} shows that it is $\uoLt_E$-convergent to zero. On the other hand, \cite[Proposition~3.1]{gao_troitsky_xanthos:2017} shows that uo-convergence of a sequence in $\Ell_1([0,1])$ is the same as almost everywhere convergence. Hence the sequence is not uo-convergent to zero.

Still, something can be salvaged in the general case. As a motivating example, suppose that $(X,\Sigma,\mu)$ is a measure space. It is well known that a sequence in $\Ell_0(X,\Sigma,\mu)$ that converges (globally) in measure has a subsequence that converges to the same limit almost everywhere; see \cite[Theorem 2.30]{folland_REAL_ANALYSIS_SECOND_EDITION:1999}, for example. When $\mu$ is finite, then, in view of  \cref{res:tau_m_is_convergence_in_measure} and \cite[Proposition~3.1]{gao_troitsky_xanthos:2017}, this can be restated as saying that a $\uoLt_E$-convergent sequence in $\Ell_0(X,\Sigma,\mu)$ has a subsequence that is uo-convergent to the same limit. We shall now extend this formulation of the result to a more general context of nets and Hausdorff \uoLtops\ on vector lattices; see \cref{res:tau_m_to_sub_uo}, below. In \cref{res:convergence_in_measure_and_almost_everywhere}, below, we shall then obtain a stronger version of the motivating result for convergence in measure and convergence almost everywhere, as a specialisation of the general result.

We start with some preparations that appear to have some independent interest.

\begin{proposition}\label{res:prep_for_tau_m_to_sub_uo_1}
	Let $E$ be a vector lattice with the countable sup property such that $\ocdual{E}$ separates the points of $E$. Take $e\in \posE$, and let $I_e$ denote the ideal that is generated in $E$ by $e$.  Then $\ocdual{(I_e)}$ separates the points of $I_e$. In fact,  there even exists a $\varphi\in \ocdual{(I_e)}$ that is strictly positive on $I_e$.
\end{proposition}

\begin{proof}
		It is immediate from \cref{res:veksler} that $\ocdual{(I_e)}$ separates the points of $I_e$. It follows from \cref{res:polars} that the ideal of $\ocdual{(I_e)}$ that is generated by a strictly positive $\varphi$ in $\ocdual{(I_e)}$ would  already separate the points of $E$.	We turn to the existence of such a strictly positive $\varphi\in\ocdual{(I_e)}$.
	
	Suppose first that $E$ is Dedekind complete.
	For $\psi\in\pos{\left(\ocdual{E}\right)}$, we let
		\[
		N_\psi\coloneqq\{x\in E: \psi(\abs{x})=0\}
		\]
		denote its null ideal, and we let
		\[
		C_\psi\coloneqq \text{N}_\psi^\text{d}
		\]
	denote its carrier. Since $\psi$ is order continuous, $N_\psi$ is a band in $E$.
	
	Let $B_0$ be the band that is generated by the subset $\{C_\psi: \psi\in \pos{\left(\ocdual{E}\right)} \}$ of $E$. Then
	\[
	 B_0^\text{d}=\bigcap_{\psi\in\pos{\left(\ocdual{E}\right)}}C_\psi^{\text{d}}=\bigcap_{\psi\in\pos{\left(\ocdual{E}\right)}}N_\psi^{\text{dd}}=\bigcap_{\psi\in\pos{\left(\ocdual{E}\right)}}N_\psi=\{0\},
	\]
	where in the final step we have used \cref{res:polars} and the fact that $\ocdual{E}$ separates the points of $E$. We thus see that $B_0=E$.
	
	For $\psi\in\pos{\left(\ocdual{E}\right)}$, let $P_{C_\psi}$ denote the band projection from $E$ onto $C_\psi$. When $\psi_1,\psi_2\in \pos{\left(\ocdual{E}\right)}$ and $\psi_1\leq\psi_2$, then $C_{\psi_1}\subseteq C_{\psi_2}$ which, by \cite[Theorem~1.46]{aliprantis_burkinshaw_POSITIVE_OPERATORS_SPRINGER_REPRINT:2006}, is equivalent to $P_{C_{\psi_1}}\leq P_{C_{\psi_2}}$. Therefore, the net $\{P_{C_\psi}:\psi \in \pos{\left(\ocdual{E}\right)}\}$ in $\regops(E)$ is increasing. Set
	\[
	P\coloneqq \sup\, \{P_{C_\psi}:\psi \in \pos{\left(\ocdual{E}\right)}\},
	\]
	where the supremum is in $\regops(E)$. From 	\cite[Theorem~30.5]{luxemburg_zaanen_RIESZ_SPACES_VOLUME_I:1971} we know that $P$ is a band projection with $B_0$ as its range space. Since $B_0=E$, it follows that $P=I$. This implies that $\{P_{C_\psi}e:\psi \in \pos{\left(\ocdual{E}\right)}\}\uparrow e$, and it follows from the fact that $E$ has the countable sup property that there exists a sequence $(\psi_n)_{n=1}^\infty$ in  $\pos{\left(\ocdual{E}\right)}$ such that $P_{C_{\psi_n}}e\uparrow e$ in $E$.
	
	Consider the ideal $I_e$ of $E$. Since $E$ is Dedekind complete it is uniformly complete, so that $I_e$ is a Banach lattice when supplied with its order unit norm $\norm{\,\cdot\,}_e$. Its order dual $\odual{I}_e$ coincides with its norm dual $E^\ast$ and is then a Banach lattice. Choose strictly positive real numbers $\alpha_1,\alpha_2,\ldots$ such that $\sum_{n=1}^{\infty} \alpha_n\norm{\psi_n|_{I_e}}<\infty$, and define $\varphi\in \odual{I}_e$ by setting
	\[
	\varphi\coloneqq \sum_{n=1}^{\infty} \alpha_n \psi_n|_{I_e}.
	\]
	Since $I_e$, being an ideal of $E$, is a regular vector sublattice of $E$, each $\psi_n|_{I_e}$ is order continuous. On observing that, being a band, $\odual{(I_e)}_{\oc}$ is an order closed and, therefore, norm closed subset of the Banach lattice $E^\ast$, we see that $\varphi$ is order continuous on $I_e$. Obviously, $\varphi$ is positive.
	
	Suppose that $x\in I_e$ is positive and that $\varphi(x)=0$. Then $\psi_n(x)=0$ for all $n\geq 1$. That is, $x\in N_{\psi_n}$ for all $n\geq 1$, so that $P_{C_{\psi_n}}x=0$ for all $n\geq 1$.
	
	Take $\lambda\geq 0$ such that $0\leq x\leq\lambda e$. Using \cite[Theorem~2.49, Theorem~2.44, and Definition~2.41]{aliprantis_burkinshaw_POSITIVE_OPERATORS_SPRINGER_REPRINT:2006}, we see that there exists an order continuous operator $T$ on $E$ that commutes with all band projections on $E$ and is such that $T(\lambda e)=x$. Since $P_{C_{\Psi_n}}(\lambda e)\uparrow \lambda e$ in $E$, we have $TP_{C_{\Psi_n}}(\lambda e)\uparrow T(\lambda e)=x$ in $E$. On the other hand, we know that $TP_{C_{\Psi_n}}(\lambda e)=P_{C_{\Psi_n}}T(\lambda e)=P_{C_{\Psi_n}}x=0$ for all $n$. We conclude that $x=0$. Hence $\varphi$ is strictly positive on $I_e$. This completes the proof when $E$ is Dedekind complete.
	
	For general $E$, we note that its Dedekind completion $E^\delta$ also has the countable sup property; see \cite[Theorem~32.9 ]{luxemburg_zaanen_RIESZ_SPACES_VOLUME_I:1971}. Furthermore, \cref{res:veksler} shows that $\ocdual{\left(E^\delta\right)}$ separates the points of $E^\delta$. Let $I_{e,\delta}$ denote the ideal that is generated by $e$ in $E^\delta$. By what has been established above, there exists a $\varphi_\delta\in \ocdual{\left(I_{e,\delta}\right)}$ that is strictly positive on $I_{e, \delta}$. Hence its restriction $\varphi_\delta|_{I_e}$ to $I_e$ is strictly positive on $I_e$. This restriction is also order continuous on $I_e$. To see this, suppose that $\net$ is a net in $I_e$ and that $x_\alpha\oconv 0$ in $I_e$. Since $I_e$, being an ideal of $E$, is a regular vector sublattice of $E$, and since $E$, being order dense in $E^\delta$, is a regular vector sublattice of $E^\delta$, $I_e$ is a regular vector sublattice of $E^\delta$. Thus $x_\alpha\oconv 0$ in $E^\delta$. There exists an $\alpha_0\in\is$ such that the tail $(x_\alpha)_{\alpha\in\is, \alpha\geq\alpha_0}$ is order bounded in $I_e$. Since this tail is then evidently also order bounded in $I_{e,\delta}$, \cref{res:local-global_for_o-convergence_and_uo-convergence} shows that $x_\alpha\oconv 0$ in $I_{e,\delta}$ for $\alpha\geq\alpha_0$. Then $\varphi_\delta|_{I_e}(x_\alpha)\conv{}0$ for $\alpha\geq\alpha_0$ by the order continuity of $\varphi$ on $I_{e,\delta}$. Consequently, $\varphi_\delta|_{I_e}(x_\alpha)\conv{}0$, as required.
\end{proof}

Suppose that a vector lattice $E$ has an order unit $e$ and that $\net$ is a net in $E$. According to \cite[Corollary~3.5]{gao_troitsky_xanthos:2017},  the fact that $\abs{x_\alpha}\wedge e\oconv 0$  is already enough to imply that  $x_\alpha\uoconv 0$. This is a special case of the following.
	
\begin{proposition}\label{res:prep_for_tau_m_to_sub_uo_3}
	Let $E$ be a vector lattice, let $S$ be a non-empty subset of $E$, and let $B_S$ denote the band that is generated by $S$ in $E$. Suppose that $\net$ is a net in $B_S$ such that $\abs{x_\alpha}\wedge \abs{y}\oconv 0$ in $E$ for all $y\in S$. Then $x_\alpha\uoconv 0$ in $E$.
\end{proposition}

\begin{proof} Since $B_S$ is a regular sublattice of $E$, part~\ref{part:order_convergence_and_regular_sublattices} of \cref{res:local-global_for_o-convergence_and_uo-convergence} shows that $\abs{x_\alpha}\wedge \abs{y}\oconv 0$ in $B_S$ for all $y\in S$. From \cite[Lemma~2.2]{li_chen:2018} it then follows that $x_\alpha\uoconv 0$ in $B_S$. Part~\ref{part:unbounded_order_convergence_and_regular_sublattices} of \cref{res:local-global_for_o-convergence_and_uo-convergence} then implies that $x_\alpha\uoconv 0$ in $E$.
\end{proof}

\begin{lemma}\label{res:countable_sup_property_and_generated_ideals}
Let $E$ be a vector lattice, and let $F$ be an order dense vector sublattice of $E$. Then $F$ has the countable sup property if and only if the ideal $I_F$ of $E$ that is generated by $F$ has the countable sup property.
\end{lemma}

\begin{proof}
Since $F$ is an order dense and majorising vector sublattice of $I_F$, the Dedekind completions of $F$ and $I_F$ are isomorphic. The proof is then completed by using that the countable sup property of a vector lattice and that of its Dedekind completion are equivalent;  see \cite[Theorem~32.9]{luxemburg_zaanen_RIESZ_SPACES_VOLUME_I:1971}.
\end{proof}

\begin{proposition}\label{res:order_separability_equivalences}
	Let $E$ be a vector lattice, and let $F$ be an order dense vector sublattice of $E$. The following are equivalent:
	\begin{enumerate}
		\item $E$ has the countable sup property;\label{part:order_separability_equivalences_1}
		\item $F$ has the countable sup property and $F$ is super order dense in $E$;\label{part:order_separability_equivalences_2}
	\end{enumerate}
\end{proposition}

\begin{proof}
	Suppose that $E$ has the countable sup property. Let the net $\net$ in $\pos{F}$ and $x$ in $\pos{F}$ be such that $x_\alpha\uparrow x$ in $F$. Since $F$, being order dense in $E$, is a regular vector sublattice of $E$, we also have that $x_\alpha\uparrow x$ in $E$. By hypothesis, there exists a sequence of indices $\seq{\alpha_n}{n}$ in $\is$ such that $x_{\alpha_n}\uparrow x$ in $E$. Then also $x_{\alpha_n}\uparrow x$ in $F$. Hence $F$ has the countable sup property. Take an $x\in F^+$. Since $F$ is order dense in $E$, \cite[Theorem~1.34]{aliprantis_burkinshaw_POSITIVE_OPERATORS_SPRINGER_REPRINT:2006} shows that $S\coloneqq\{y\in F:0\leq y\leq x\}\uparrow x$ in $E$. The fact that $E$ has the countable sup property then yields a sequence $\seq{x_n}{n}\subseteq S\subseteq F$ such that $x_n\uparrow x$ in $E$. Hence $F$ is super order dense in $E$.
	
	Suppose that $F$ has the countable sup property and that $F$ is super order dense in $E$. According to \cref{res:countable_sup_property_and_generated_ideals}, the ideal $I_F$ in $E$ that is generated by $F$ also has the countable sup property. Since $I_F\supseteq F$ is evidently super order dense in $E$, it follows from \cite[Theorem~29.4]{luxemburg_zaanen_RIESZ_SPACES_VOLUME_I:1971} that $E$ has the countable sup property.
\end{proof}	

\begin{remark}\quad\begin{enumerate}
		\item The assumption in part~\ref{part:order_separability_equivalences_2} that $F$ be super order dense in $E$ cannot be relaxed to requiring it to be merely order dense in $E$, not even when $F$ is an  order dense ideal of $E$ rather than an order dense vector sublattice; see \cite[Example~4.3]{kandic_vavpetic:2018}.
		\item When $E$ is a vector lattice and $F$ is a vector sublattice of $E$ that has the countable sup property, then the fact that $E$ has the countable sup property is equivalent to the super order density of the Dedekind completion of $F$ in the Dedekind completion of $E$. We refer to \cite[Theorem~4.5]{kandic_vavpetic:2018} for this result in the same spirit as \cref{res:order_separability_equivalences}.
	\end{enumerate}
\end{remark}

All preparations have now been made for the proof of the core result of this section.

\begin{theorem}\label{res:tau_m_to_sub_uo}
	
	Let $E$ be a vector lattice with the countable sup property, and suppose that $E$ has an order dense ideal $F$ such that $\ocdual{F}$ separates the points of $F$. Let $G$ be a regular vector sublattice of $E$. Then $G$ admits a \uppars{necessarily unique} Hausdorff \uoLtop\ $\uoLt_G$.
	
	Let $\net$ be a net in $G$ and suppose that $x_\alpha\conv{\uoLt_G} x$ for some $x\in G$. Take a sequence $(\alpha^\prime_n)_{n=1}^\infty$ of indices in $\is$. Then there exists an increasing sequence $\alpha_1^\prime=\alpha_1\leq\alpha_2\leq\dotsb$ of indices in $\is$ such that $\alpha_n\geq\alpha^\prime_n$ for all $n\geq 1$ and  $x_{\alpha_n}\uoconv x$ in $G$. In particular, when a sequence  $(x_n)_{n=1}^\infty $ in $G$ and $x\in G$ are such that $x_n\conv{\uoLt_G}x$ in $G$, then there exists a subsequence $(x_{n_k})_{k=1}^\infty$ of $(x_n)_{n=1}^\infty $ such that $x_{n_k}\uoconv x$ in $G$.
\end{theorem}

\begin{proof}
	
	In view of \cref{res:going_down_one} and \cref{res:local-global_for_o-convergence_and_uo-convergence}, we may (and shall) suppose that $G=E$.
	
	We know from \cref{res:uo-lebesgue_topology_originating_from_separating_order_continuous_dual_of_order_dense_ideal} that $E$ admits a (necessarily unique) Hausdorff \uoLtop\ $\uoLt_E$.
		The statement on subsequences is clear from the statement on nets, so we need only establish the existence of the $\alpha_n$ for $n\geq 1$. We may suppose that $x=0$.

	Suppose first that $E$ is Dedekind complete.

	For $y\in \posF$, we let $I_y\subseteq F$ denote the ideal that is generated by $y$ in $E$. By \cref{res:order_separability_equivalences}, $F$ inherits the countable sup property from $E$. Hence  \cref{res:prep_for_tau_m_to_sub_uo_1} applies to the vector lattice $F$.  We then see that $(I_y)^\sim_\oc$ separates the points of $I_y$ and that there even exists a strictly positive order continuous linear functional on $I_y$. We choose and fix such a strictly positive $\varphi_y\in(I_y)^\sim_\oc$ for each $y\in \posF$.
	From \cref{res:essential} we know that
	\begin{equation}\label{eq:converges_to_zero}
	\varphi_y(\abs{x_\alpha}\wedge y)\to 0
	\end{equation}
	for all $y\in \posF$.

	Set $\alpha_1\coloneqq\alpha_1^\prime$. Since $F$ is super order dense in $E$ by \cref{res:order_separability_equivalences}, we can choose a sequence $\{y_m^1\}_{m=1}^\infty$ in $\posF$ such that $0\leq y_m^1\uparrow_m \abs{x_{\alpha_1}}$.
	
	For $n\geq 2$,  we shall now inductively construct an index $\alpha_n\in\is$  and a sequence $\{y_m^n\}_{m=1}^\infty$ in $\posF$  such that, for all $n\geq 2$:
	\begin{enumerate_alpha}
		\item $\alpha_n\geq\alpha_n^\prime$;
		\item $\alpha_{n}\geq\alpha_{n-1}$;
		\item $\varphi_{y_m^i}\bigl(\abs{x_{\alpha_n}}\wedge y_m^i\bigr)<2^{-n}$ for $i=1,2,\ldots,n-1$ and $m=1,2,\dotsc,n$;
		\item $0\leq y_m^n\uparrow_m \abs{x_{\alpha_n}}$ in $E$.
	\end{enumerate_alpha}
	
	We start with $n=2$. The elements $y_m^1$ of $\posF$ are already known for all $m\geq 1$, and $\varphi_{y_m^1}\bigl(\abs{x_\alpha}\wedge y_m^1\bigr)\conv{}0$ for all $m\geq 1$ by \cref{eq:converges_to_zero}. Therefore, we can choose an $\alpha_2\in\is$ such that  $\varphi_{y_m^1}\bigl(\abs{x_{\alpha_2}}\wedge y_m^1\bigr)<2^{-2}$ for $m=1,2$. We can arrange that also $\alpha_2\geq\alpha_2^\prime$ and $\alpha_2\geq\alpha_1$. Finally, we choose a sequence $(y_m^2)_{m=1}^\infty$ in $F$ such that $0\leq y_m^2\uparrow_m \abs{x_{\alpha_2}}$. This completes the construction for $n=2$.
	
	Suppose that $n\geq 2$ and that we have already constructed $\alpha_2,\ldots,\alpha_n\in\is$ and sequences $(y_m^1)_{m=1}^\infty,\dotsc, (y_m^n)_{m=1}^\infty$ in $\posF$ satisfying the four requirements above. The elements $y_m^i$ of $\posF$ are already known for all $i=1,2,\ldots,n$ and $m\geq 1$, and $\varphi_{y_m^i}\bigl(\abs{x_\alpha}\wedge y_m^i\bigr)\conv{}0$ for all such $i$ and $m$ by \cref{eq:converges_to_zero}. Therefore, we can choose $\alpha_{n+1}\in\is$ such that $\varphi_{y_m^i}\bigl(\abs{x_{\alpha_{n+1}}}\wedge y_m^i\bigr)<2^{-(n+1)}$ for all $i=1,2,\ldots,n$ and $m=1,2,\dotsc,n+1$. We can arrange that also $\alpha_{n+1}\geq\alpha_{n+1}^\prime$ and $\alpha_{n+1}\geq\alpha_n$. Finally, we choose a sequence $(y_m^{n+1})_{m=1}^\infty$ in $\posF$ such that $0\leq y_m^{n+1}\uparrow_m \abs{x_{\alpha_{n+1}}}$ in $E$. This completes the construction for $n+1$.

	Fix $i,m\geq 1$. Since $0\leq \abs{x_{\alpha_j}}\wedge y_m^i\leq y_m^i$ for all $j\geq 1$, we can define elements $z_n^{j,m}$ of $I_{y_m ^i}$ for $n\geq 1$ by setting  $z_{n}^{i,m}\coloneqq\bigvee_{j=n}^{\infty}\bigl(\abs{x_{\alpha_j}}\wedge y_m^i\bigr)$. Here the  supremum is in the ideal $I_{y_m ^i}$ in $E$ (which, although this is immaterial, happens to coincide with the supremum in $E$). It is clear that $z_n\geq 0$ for $n\geq 1$ and that $z_{n}^{i,m}\downarrow_n$; we shall show that $z_{n}^{i,m}\downarrow_n 0$ in $I_{y_m^i}$. For this, we start by noting that
	the inequality in (c) shows that $\varphi_{y_m^i}\bigl(\abs{x_{\alpha_j}}\wedge y_m^i\bigr)<2^{-j}$ for all $j\geq \max(i+1,m)$. Therefore, for all $n\geq \max(i+1,m)$, we can use the order continuity of $\varphi_{y_m^i}$ on $I_{y_m^i}$ to see that
	\begin{align*}
	0&\leq \varphi_{y_m^i}(z_{n}^{i,m})
	\\&=\varphi_{y_m^i}\left(\bigvee_{j=n}^{\infty}\bigl(\abs{x_{\alpha_j}}\wedge y_i^m\bigr)\right)
	\\ &=\varphi_{y_m^i}\left(\sup_{k\geq n}\left(\bigvee_{j=n}^{k}\bigl(\abs{x_{\alpha_j}}\wedge y_i^m\bigr)\right)\right)
	\\ &=\lim_{\overset{k\to\infty}{k\geq n}}\varphi_{y_m^i}\left(\bigvee_{j=n}^{k}\bigl(\abs{x_{\alpha_j}}\wedge y_i^m\bigr)\right)\\
	\\ &\leq \limsup_{\overset{k\to\infty}{k\geq n}}\varphi_{y_m^i}\left(\sum_{j=n}^{k}\bigl(\abs{x_{\alpha_j}}\wedge y_i^m\bigr)\right)
	\\ &\leq \limsup_{\overset{k\to\infty}{k\geq n}}\sum_{j=n}^k 2^{-j}
	\\&\leq 2^{-n+1}.
	\end{align*}
	
	We see from this that for the infimum $\inf_{n\geq 1}z_n^{i,m}$ in $I_{y_m^i}$ (which, although again immaterial, happens to coincide with the infimum in $E$) we have $0\leq \varphi_{y_m^i}\left(\inf_{n\geq 1} z_n^{i,m}\right)\leq 2^{-n+1}$ for all $n\geq \max(i+1,m)$. Hence $\varphi_{y_m^i}\left(\inf_{n\geq 1} z_n^{i,m}\right)=0$. Since $\varphi_{y_m^i}$ is strictly positive on $I_{y_m^i}$, this implies that $\inf_{n\geq 1}z_n^{i,m}=0$ in $I_{y_m^i}$, as we wanted to show.
	
	The inequalities $0\leq \abs{x_{\alpha_n}}\wedge y_m^i\leq z_n^{i,m}$ for all $n\geq 1$ now show that $\abs{x_{\alpha_n}}\wedge y_m^i\oconv 0$ in $I_{y_m^i}$ as $n\to\infty$, and then also  $\abs{x_{\alpha_n}}\wedge y_m^i\oconv 0$ in $E$ as $n\to\infty$.
	
	We have now shown that, for all $i,m\geq 1$, $\abs{x_{\alpha_n}}\wedge y_m^i\oconv 0$ in $E$ as $n\to\infty$.
	
	Let $B$ denote the band that is generated by $\{y_m^i:i,m\geq 1\}$ in $E$. In view of (d) above, it is clear that the sequence $(x_{\alpha_n})_{n=1}^\infty$ is a sequence in $B$. We can now conclude from \cref{res:prep_for_tau_m_to_sub_uo_3} that $x_{\alpha_n}\uoconv 0$ in $E$.  This concludes the proof when $E$ is Dedekind complete.
	
	For a general vector lattice $E$, we pass to the Dedekind completion $E^\delta$ of $E$. By \cite[Theorem~32.9]{luxemburg_zaanen_RIESZ_SPACES_VOLUME_I:1971}, $E^\delta$ also has the countable sup property. We let $F^\delta$ denote the ideal that is generated in $E^\delta$ by $F$. Then $F$ is obviously majorising in $F^\delta$. Since $F$ is order dense in $E$ and $E$ is order dense in $E^\delta$, $F$ is order dense in $E^\delta$ and then also in $F^\delta$. We see from this that, as the notation  already suggests, $F^\delta$ is the Dedekind completion of $F$, but what we actually need is that, by \cref{res:veksler},  $\ocdual{\left(F^\delta\right)}$ separates the points of $F^\delta$. The fact that $F$ is order dense in $E^\delta$ implies that $F^\delta\supseteq F$ is order dense in $E^\delta$. Hence $E^\delta$ also admits a (necessarily) unique Hausdorff \oLtop\ $\uoLt_{E^\delta}$. Moreover, \cref{res:going_down_one} shows that $x_\alpha\conv{\uoLt_{E^\delta}}0$ in $E^\delta$. 	
	By what has been established for the Dedekind complete case, there exist indices $\alpha_n$ as specified such that $x_{\alpha_n}\uoconv 0$ in $E^\delta$.  By \cref{res:local-global_for_o-convergence_and_uo-convergence}, $x_{\alpha_n}\uoconv 0$ in $E$.
\end{proof}

For comparison, we include the following; see \cite[Theorem~4.19]{aliprantis_burkinshaw_LOCALLY_SOLID_RIESZ_SPACES_WITH_APPLICATIONS_TO_ECONOMICS_SECOND_EDITION:2003}. We recall that a topology on a vector lattice $E$ is a \emph{Fatou topology} when it is a (not necessarily Hausdorff) locally solid linear topology on $E$ that has a base of neighbourhoods of zero consisting of solid and order closed sets. A Lebesgue topology is a Fatou topology; see \cite[Lemma~4.1]{aliprantis_burkinshaw_LOCALLY_SOLID_RIESZ_SPACES_WITH_APPLICATIONS_TO_ECONOMICS_SECOND_EDITION:2003}, for example.

\begin{theorem}\label{res:other_subsequence_result}
	Let $E$ be a vector lattice with the countable sup property that is supplied with a Hausdorff locally solid  linear topology $\tau$ with the Fatou property. Suppose that $\net$ is an order bounded net in $E$ and that $x_\alpha\conv{\tau} x$ for some $x\in E$. Then there exist indices $\alpha_1\leq\alpha_2\leq\dotsb$ in $\is$ such that $x_{\alpha_n}\oconv x$.
\end{theorem}

The hypotheses in \cref{res:other_subsequence_result} on the topology on the vector lattice are weaker than those in   \cref{res:tau_m_to_sub_uo}, and its conclusion is stronger. The big difference is, however, that the net in \cref{res:other_subsequence_result} is supposed to be order bounded, whereas there is no such restriction in \cref{res:tau_m_to_sub_uo}.

\cref{res:other_subsequence_result} also holds when, instead of requiring $E$ to have the countable sup property, it is required that there exist an at most countably infinite subset of $E$ such that the band  that it generates equals the carrier of $\tau$; see \cite[Theorem~6.7]{kandic_taylor:2018}. We refer to \cite[Definition~4.15]{aliprantis_burkinshaw_LOCALLY_SOLID_RIESZ_SPACES_WITH_APPLICATIONS_TO_ECONOMICS_SECOND_EDITION:2003} for the definition of the carrier of a (not necessarily Hausdorff) locally solid topology on a vector lattice.

\begin{remark}
	The hypothesis in \cref{res:tau_m_to_sub_uo} that $E$ have the countable sup property cannot be relaxed to merely requiring that $F$ have this property.  As a counter-example, consider the situation where $F$ is a Banach lattice with an order continuous norm that is an order dense ideal of a vector lattice $E$. Then $\ocdual{F}=F^\ast$ separates the points of $F$, and it is easy to see that $F$ has the countable sup property; the latter also follows from a more general result in  \cite[Theorem~4.26]{aliprantis_burkinshaw_LOCALLY_SOLID_RIESZ_SPACES_WITH_APPLICATIONS_TO_ECONOMICS_SECOND_EDITION:2003}.
	Since the norm topology on $F$ is a Hausdorff \oLtop\ on $F$, $E$ has a (necessarily unique) Hausdorff \uoLtop\ $\uoLt_E$. It is the topology of un-convergence with respect to $F$. It is possible to find such $F$ and $E$, and a sequence in $E$ that is $\uoLt_E$-convergent to zero in $E$, yet has no subsequence that is uo-convergent to zero in $E$;  see \cite[Example 9.6]{kandic_li_troitsky:2018}.
\end{remark}

We have the following consequence of \cref{res:tau_m_is_convergence_in_measure} and \cref{res:tau_m_to_sub_uo}.

\begin{theorem}\label{res:convergence_in_measure_and_almost_everywhere}
	Let $(X,\Sigma,\mu)$ be a measure space where $\mu$ is $\sigma$-finite. Suppose that $(f_\alpha)_{\alpha\in\is}$ is a net in $\Ell_0(X,\Lambda,\mu)$ such that $f_\alpha\conv{\mu^\ast}0$. Take a sequence $(\alpha^\prime_n)_{n=1}^\infty$ of indices in $\is$. Then there exists an increasing sequence $\alpha_1^\prime=\alpha_1\leq\alpha_2\leq\dotsb$ of indices in $\is$ such that $\alpha_n\geq\alpha^\prime_n$ for all $n\geq 1$ and  $f_{\alpha_n}\conv{} 0$ almost everywhere. In particular, when a sequence  $(f_n)_{n=1}^\infty$ is a sequence in $\Ell_0(X,\Lambda,\mu)$ and $f_n\conv{\mu^\ast}0$, then there exists a subsequence $(f_{n_k})_{k=1}^\infty$ of $(f_n)_{n=1}^\infty $ such that $f_{n_k}\conv{}0$ almost everywhere.
\end{theorem}

\begin{proof}
	It is known that $\Ell_0(X,\Sigma,\mu)$ has the countable sup property for every $\sigma$-finite measure $\mu$; see \cite[Theorem~7.73]{aliprantis_burkinshaw_LOCALLY_SOLID_RIESZ_SPACES_WITH_APPLICATIONS_TO_ECONOMICS_SECOND_EDITION:2003} or \cite[Lemma~2.6.1]{meyer-nieberg_BANACH_LATTICES:1991}, for example.
	
	The combination of \cref{res:tau_m_is_convergence_in_measure} and \cref{res:tau_m_to_sub_uo} yields a sequence of indices $\alpha_n$ as specified such that $f_{\alpha_n}\uoconv 0$. Since, for a general measure $\mu$, uo-convergence of a sequence in $\Ell_0(X,\Sigma,\mu)$ is equivalent to its convergence almost everywhere (see \cite[Proposition~3.1]{gao_troitsky_xanthos:2017}), the proof is complete.
\end{proof}

\begin{remark}\quad
	\begin{enumerate}
		\item In view of its proof, the natural condition on $\mu$ in \cref{res:convergence_in_measure_and_almost_everywhere} is that $\mu$ be semi-finite and have the countable sup property. It is known, however, that this is equivalent to requiring that $\mu$ be $\sigma$-finite; see \cite[Proposition~6.5]{kandic_taylor:2018}.
		\item For every measure $\mu$, a sequence in $\Ell_0(X,\Lambda,\mu)$ that converges (globally) in measure has a subsequence that converges almost everywhere to the same limit; see \cite[Theorem 2.30]{folland_REAL_ANALYSIS_SECOND_EDITION:1999}, for example. \cref{res:convergence_in_measure_and_almost_everywhere} does not imply this result for arbitrary measures, but once the measure is known to be $\sigma$-finite, it \emph{does} produce the desired subsequence, and it even does so under the weaker hypothesis of convergence in measure on subsets of finite measure.
	\item Even for finite measures, we are not aware of an existing result that, as in \cref{res:convergence_in_measure_and_almost_everywhere},  is concerned with \emph{nets} that converge in measure.
\end{enumerate}

\end{remark}

We conclude this section by extending another classical result from measure theory to the context of Hausdorff \uoLtops\ and uo-convergence. Suppose that $(X,\Sigma.\mu)$ is a measure space, where $\mu$ is $\sigma$-finite. Then a sequence in $\Ell_0(X,\Sigma,\mu)$ is convergent in measure on subsets of finite measure if and only if every subsequence has a further subsequence that converges to the same limit almost everywhere; see \cite[Exercise~18.14 on p.~132]{zaanen_INTEGRATION:1967}. This is a special case of the following.

\begin{theorem}\label{res:tau_m_to_sub_uo_iff}
	
		Let $E$ be a vector lattice with the countable sup property, and suppose that $E$ has an order dense ideal $F$ such that $\ocdual{F}$ separates the points of $F$. Let $G$ be a regular sublattice of $E$. Then $G$ admits a \uppars{necessarily unique} Hausdorff \uoLtop\ $\uoLt_G$.
		For a sequence $\seq{x_n}{n}\subseteq G$, $x_n \conv{\uoLt_G} 0$ in $G$ if and only if every subsequence $\seq{x_{n_k}}{k}$ of $\seq{x_n}{n}$ has a further subsequence $\seq{x_{n_{k_i}}}{i}$ such that $x_{n_{k_i}}\uoconv 0$ in $G$.
\end{theorem}

\begin{proof}
	
	In view of \cref{res:going_down_one} and \cref{res:local-global_for_o-convergence_and_uo-convergence}, we may (and shall) suppose that $G=E$.
	
	The forward implication is clear from \cref{res:tau_m_to_sub_uo}.  We now show the converse. When it fails that $x_n\uoLtconv 0$ in $E$, then \cref{res:uo-lebesgue_topology_originating_from_separating_order_continuous_dual_of_order_dense_ideal} shows that there exists an $\varphi\in\ocdual{F}$, an $y\in F$, a subsequence $\seq{x_{n_k}}{k}$ of $\seq{x_n}{n}$ and an $\varepsilon>0$ such that $\abs{\varphi}(\abs{x_{n_k}}\wedge\abs{y})> \varepsilon$ for all $k$. It is then clear from the order continuity of $\varphi$ that it is impossible to find a further subsequence $\seq{x_{n_{k_i}}}{i}$ of $\seq{x_{n_k}}{k}$  such that $x_{n_{k_i}}\uoconv 0$ in $E$.
\end{proof}

As another special case of \cref{res:tau_m_to_sub_uo_iff}, we see that a sequence in a Banach lattice with an order continuous norm is un-convergent to zero if and only if every subsequence has a further subsequence that is uo-convergent to zero; we recall that a Banach lattice with an order continuous norm has the countable sup property by \cite[Theorem~4.26]{aliprantis_burkinshaw_LOCALLY_SOLID_RIESZ_SPACES_WITH_APPLICATIONS_TO_ECONOMICS_SECOND_EDITION:2003}, for example. We have thus retrieved \cite[Theorem 4.4]{deng_o_brien_troitsky:2017}.

\section{Topological aspects of (unbounded) order convergence}\label{sec:topological aspects of uo-convergence}

\noindent In this section, we consider topological issues that are related to (sequential) order convergence and to (sequential) unbounded order convergence, with an emphasis on the latter. \cref{res:tau_m_to_sub_uo} will be seen to be an important tool.

Let $E$ be a vector lattice, and let $A\subseteq E$. We define the \emph{\oadhtext\ of $A$} as the set of all order limits of nets in $A$, and denote it by $\oadh{A}$. The \emph{\soadhtext\ of $A$} is the set of all order limits of sequences in $A$; it is denoted by $\soadh{A}$. \footnote{In  \cite[p.~82]{luxemburg_zaanen_RIESZ_SPACES_VOLUME_I:1971}, our \soadhtext\ is called the pseudo order closure. In \cite{gao_leung:2018}, our \oadhtext\ of a subset $A$ is called the order closure of $A$, and it is denoted by $\overline{A}^\o$. These two terminologies, as well as the notation $\overline{A}^\o$, could suggest that taking the (pseudo) order closure is a (sequential) closure operation for a topology. Since this is hardly ever the case, we prefer a terminology and notation that avoid this possible confusion. It is inspired by \cite[Definition~1.3.1]{beattie_butzmann_CONVERGENCE_STRUCTURES_AND_APPLICATIONS_TO_FUNCTIONAL_ANALYSIS:2002}.} The \emph{\uoadhtext}  $\uoadh{A}$ and the \emph{\suoadhtext} $\suoadh{A}$ of $A$ are similarly defined. The subset $A$ is \emph{o-closed} when $\oadh{A}=A$.\footnote{This definition is consistent with that in \cite{gao_leung:2018}.} The collection of all o-closed subsets of $E$ is easily seen to be the collection of closed sets of a topology that is called the \emph{o-topology on $E$}. The closure of a subset $A$ in the o-topology is denoted by $\oclos{A}$.\footnote{There is no notation for the closure operation in the o-topology in \cite{gao_leung:2018}.} We have $\oadh{A}\subseteq\oclos{A}$, with equality if and only if $\oadh{A}$ is o-closed. Likewise, there are $\sigma$o-closed subsets and a $\sigma$o-topology, uo-closed subsets and a uo-topology, and $\sigma$uo-closed subsets and a $\sigma$uo-topology, with similar notations and statements about inclusions and equalities of sets. Evidently, a uo-closed subset is o-closed, and a $\sigma$uo-closed subset is $\sigma$o-closed.

Order convergence in a vector lattice $E$ is hardly ever topological; according to \cite[Theorem~2.2]{dabboorasad_emelyanov_marabeh:2020} or \cite[Theorem~18.36]{taylor_THESIS:2018}, this is the case if and only if $E$ is finite-dimensional. It is not even true that the set map $A\mapsto\oadh{A}$ is always idempotent, i.e., that the \oadhtext\ of a set is always o-closed. It is known, for example, that in every $\sigma$-order complete Banach lattice that does \emph{not} have an order continuous norm, there even exists a vector sublattice such that its \oadhtext\ is not order closed; see \cite[Theorem~2.7]{gao_leung:2018}.

We know from \cref{res:atomic} that uo-convergence in atomic vector lattices is topological. According to \cite[Theorem~6.45]{taylor_THESIS:2018}, atomic vector lattices are, in fact, the only ones for which this is the case.

It appears to be open whether the \uoadhtext\ of a subset of a vector lattice is always uo-closed. In \cite[Problem~2.5]{gao_leung:2018}, it is even asked whether the \uoadhtext\ of a vector sublattice is always o-closed, which is asking for a weaker conclusion for a much more restrictive class of subsets.

Even though the topological aspects of uo-convergence are still not well understood in general, there is a class of vector lattices where we have a reasonably complete picture. In order to formulate this, we need some more notation. For a set $X$ with a topology $\tau$ and a subset $A\subseteq X$ of $X$, we let $\sadh{\tau}{A}$ denote the $\sigma\tau$-adherence of $A$, i.e., $\sadh{\tau}{A}$ is the set consisting of all $\tau$-limits of sequences in $A$. When $\sadh{\tau}{A}=A$, $A$ is said to be $\sigma\tau$-closed. The $\sigma\tau$-closed subsets of $X$ are the closed subsets of a topology on $X$ that is called the $\sigma\tau$-topology on $X$. We let $\overline{A}^\tau$ and $\overline{A}^{\sigma\tau}$ denote the $\tau$-closure and the $\sigma\tau$-closure of a subset $A$ of $X$, respectively. Then $\sadh{\tau}{A}\subseteq\overline{A}^{\sigma\tau}$, with equality if and only if $\sadh{\tau}{A}$ is $\sigma\tau$-closed.

\begin{theorem}\label{res:seven_sets_equal}
Let $E$ be a vector lattice with the countable sup property, and suppose that $E$ has an order dense ideal $F$ such that $\ocdual{F}$ separates the points of $F$. Let $G$ be a regular vector sublattice of $E$. Then $G$ admits a \uppars{necessarily unique} Hausdorff \uoLtop\ $\uoLt_G$. For a subset $A$ of $G$, the following seven subsets of $G$ are all equal:

\begin{enumerate}
	\item $\sadh{\uoLt_G}{A}$ and $\overline{A}^{\sigma\uoLt_G}$;
	\item $\suoadh{A}$ and $\suoclos{A}$;
	\item $\uoadh{A}$ and $\uoclos{A}$;
	\item $\overline{A}^{\uoLt_G}$.
\end{enumerate}

In particular, the $\sigma\uoLt_G$-topology, the $\sigma$uo-topology, and the uo-topology on $G$ all coincide with $\uoLt_G$.
\end{theorem}

In \cref{res:seven_sets_equal}, the topological closures and ($\sigma$-)adherences are to be taken with respect to the topologies and convergences in $G$.

\begin{proof}
The existence and uniqueness of $\uoLt_G$ are clear from \cref{res:tau_m_to_sub_uo}.
Using \cref{res:tau_m_to_sub_uo} for the first inclusions, we have, for an arbitrary subset $A$ of $G$,
\[
 \overline{A}^{\uoLt_G}\subseteq\suoadh{A}\subseteq\uoadh{A}\subseteq\overline{A}^{\uoLt_G}
\]
and
\[
\sadh{\uoLt_G}{A}\subseteq \suoadh{A}\subseteq\sadh{\uoLt_G}{A}.
\]
This gives equality of $\sadh{\uoLt_G}{A}$, $\suoadh{A}$, $\uoadh{A}$, and $\overline{A}^{\uoLt_G}$. Since the set map $A\mapsto \overline{A}^{\uoLt_G}$ is idempotent, so is $A\mapsto\sadh{\uoLt_G}{A}$. Hence $\sadh{\uoLt_G}{A}$ is $\sigma\uoLt_G$-closed, so that it coincides with the $\sigma\tau$-closure $\overline{A}^{\sigma\uoLt_G}$ of $A$. A similar argument works for $\suoclos{A}$ and $\uoclos{A}$.
\end{proof}

\begin{remark}
Taking $G=E$ in \cref{res:seven_sets_equal}, the equality of $\overline{A}^{\uoLt_G}$ and $\suoadh{A}$ implies that, for a $\sigma$-finite measure $\mu$, a subset of $\Ell_0(X,\Sigma,\mu)$ is closed in the topology of convergence in measure on subsets of finite measure if and only if it contains the almost every limits of sequences in it. This is \cite[245L(b)]{fremlin_MEASURE_THEORY_VOLUME_2:2003}.
\end{remark}

In the context of \cref{res:seven_sets_equal}, it is also possible to give a necessary and sufficient condition for sequential uo-convergence to be topological; see \cref{res:sequential_uo_convergence_topological_two}, below.
The proof of the following preparatory lemma is an abstraction of the argument in \cite{ordman:1966}.

\begin{lemma}\label{res:ordman_argument}
	 Let $E$ be a vector lattice that is supplied with a topology $\tau$. Suppose that $\tau$ has the following properties:
	\begin{enumerate}
		\item for every sequence $\seq{x_n}{n}$ in $E$ and for every $x\in E$, the fact that $x_n\tauconv x$ implies that there exists a subsequence $\seq{x_{n_k}}{k}$ of $\seq{x_n}{n}$ such that $x_{n_k}\uoconv x$ as $k\to\infty$.
		\item there exists a sequence $\seq{x_n}{n}$ in $E$ and an $x\in E$ such that $x_n\tauconv x$ but $x_n\overset{\uo}{\nrightarrow} x$;
	\end{enumerate}
Then there does not exist a topology $\tau^\prime$ on $E$ such that,  for every sequence $\seq{x_n}{n}$ in $E$ and for every $x\in E$, $x_n\uoconv x$ if and only if $x_n\conv{\tau^\prime} x$.
\end{lemma}

\begin{proof}
Suppose that there were such a topology $\tau^\prime$. Take a sequence $\seq{x_n}{n}$ in $E$ and an $x\in E$ such that $x_n\tauconv x$ but $x_n\overset{\uo}{\nrightarrow} x$. Then also $x_n\overset{\tau^\prime}{\nrightarrow} x$, so that there exists a $\tau^\prime$-neighbourhood $V$ of $x$ and a subsequence $\seq{x_{n_k}}{k}$ of $\seq{x_n}{n}$ such that $x_{n_k}\not\in V$ for all $k\geq 1$. Since also $x_{n_k}\tauconv x$ as $k\to\infty$, there exists a subsequence $\seq{x_{n_{k_i}}}{i}$ of $\seq{x_{n_k}}{k}$ such that $x_{n_{k_i}}\uoconv x$ as $i\to\infty$. Hence also $x_{n_{k_i}}\conv{\tau^\prime}x$ as $i\to\infty$. But this is impossible, since the entire sequence $\seq{x_{n_{k_i}}}{i}$ stays outside $V$.
\end{proof}

The following is a direct consequence of \cref{res:ordman_argument}. The topology $\tau$ in it could be a \uoLtop, but for the result to hold it need not even be a linear topology, nor need the topology $\tau^\prime$ be.

\begin{proposition}\label{res:sequential_uo_convergence_topological_one}
	Let $E$ be a vector lattice that is supplied with a topology $\tau$. Suppose that $\tau$ has the following properties:
		\begin{enumerate}
			\item for every sequence $\seq{x_n}{n}$ in $E$ and for every $x\in E$, the fact that $x_n\uoconv x$ implies that $x_n\tauconv x$;
			\item for every sequence $\seq{x_n}{n}$ in $E$ and for every $x\in E$, the fact that $x_n\tauconv x$ implies that there exists a subsequence $\seq{x_{n_k}}{k}$ of $\seq{x_n}{n}$ such that $x_{n_k}\uoconv x$ as $k\to\infty$.
		\end{enumerate}
	Then the following are equivalent;
	\begin{enumerate}
		\item there exists a topology $\tau^\prime$ on $E$ such that, for every sequence $\seq{x_n}{n}$ in $E$ and for every $x\in E$, $x_n\uoconv x$ if and only if $x_n\conv{\tau^\prime} x$;
		\item for every sequence $\seq{x_n}{n}$ in $E$ and for every $x\in E$, the fact that $x_n\tauconv x$ implies that $x_n\uoconv x$.
	\end{enumerate}
In that case, one can take $\tau$ for $\tau^\prime$.
\end{proposition}

In the appropriate context, the combination of \cref{res:tau_m_to_sub_uo} and \cref{res:sequential_uo_convergence_topological_one} yields the following necessary and sufficient condition for sequential uo-convergence to be topological. Note that there are no assumptions at all on the topology $\tau$ in its first part.

\begin{corollary}\label{res:sequential_uo_convergence_topological_two}
Let $E$ be a vector lattice with the countable sup property, and suppose that $E$ has an order dense ideal $F$ such that $\ocdual{F}$ separates the points of $F$. Let $G$ be a regular vector sublattice of $E$. Then $G$ admits a \uppars{necessarily unique} Hausdorff \uoLtop\ $\uoLt_G$, and the following are equivalent:
\begin{enumerate}
	\item there exists a topology $\tau$ on $G$ such that, for every sequence $\seq{x_n}{n}$ in $G$ and for every $x\in G$, $x_n\uoconv x$ in $G$ if and only if $x_n\conv{\tau} x$;\label{part:uo_convergence_implies_top_convergence}
	\item for every sequence $\seq{x_n}{n}$ in $G$ and for every $x\in G$, the fact that $x_n\conv{\uoLt_G} x$ in $G$ implies that $x_n\uoconv x$ in $G$.\label{part:top_convergence_implies_uo_convergence}
\end{enumerate}
In that case, one can take $\uoLt_G$ for $\tau$.
\end{corollary}

The proof of the following result closely follows the one in \cite{ordman:1966}, where it is shown that sequential almost everywhere convergence in $L_\infty([0,1])$ is not topological.

\begin{corollary}\label{res:ordman_construction}
	 Let $(X,\Sigma,\mu)$ be a measure space, where $\mu$ is $\sigma$-finite. Suppose that there exists an $A\in\Sigma$ with the property that, for every $k\geq 1$, there exist finitely many mutually disjoint $A_{k,1},\dotsc, A_{k,N_k}\in\Sigma$ such that $0<\mu(A_{k,1}),\dotsc,\mu(A_{k,N_k})<1/k$ and $A=\bigcup_{l=1}^{N_k}A_{k,l}$.
	
Take a regular vector sublattice $G$ of $\Ell_0(X,\Sigma,\mu)$ that contains the characteristic functions $1_{A_{k,l}}$ of all sets $A_{k,l}$ for $k=1,2,\dotsc$ and $l=1,\dotsc,N_k$.  Then there does not exist a topology $\tau$ on $G$ such that, for every sequence $\seq{x_n}{n}$ in $G$ and for every $x\in G$, $x_n\uoconv x$ in $G$ if and only if $x_n\conv{\tau} x$.
\end{corollary}

\begin{proof}
 We are in the situation of \cref{res:sequential_uo_convergence_topological_two}, where $\uoLt_G$-convergence is convergence in measure on subsets of finite measure by \cref{res:tau_m_is_convergence_in_measure}, and sequential uo-convergence is almost everywhere convergence by \cite[Proposition~3.1]{gao_troitsky_xanthos:2017}. Consider the following sequence in $G$:
\[
A_{1,1},\dotsc, A_{1,N_1},A_{2,1},\dotsc, A_{2,N_2},A_{3,1},\dotsc,A_{3,N_3},\dotsc.
\]
This sequence clearly converges to zero on subsets of finite measure, but it converges nowhere to zero on the subset $A$ of strictly positive measure. Hence the property in part~\ref{part:top_convergence_implies_uo_convergence} of \cref{res:sequential_uo_convergence_topological_two} does not hold, and then neither does the property in its part~
\ref{part:uo_convergence_implies_top_convergence}.
\end{proof}	

\begin{remark}
\cref{res:ordman_construction} provides us with a large class of examples of vector lattices where sequential uo-convergence is not topological\textemdash so that uo-convergence is certainly not topological\textemdash but where, according to \cref{res:seven_sets_equal}, the set maps $A\mapsto\suoadh{A}$ and $A\mapsto\uoadh{A}$ are both still idempotent, so that $\suoadh{A}$ is $\sigma$uo-closed and $\uoadh{A}$ is uo-closed for every subset $A$ of $G$. For all $p$ such that $0\leq p\leq\infty$, the space $L_p([0,1])$ is such an example.
\end{remark}
      	
 We conclude with a strengthened version of \cite[Theorem~2.2]{gao_leung:2018}. The improvement lies in the removal of the hypothesis that $E$ be Banach lattice, and by adding eight more equal, but not obviously equal,  sets to the three equal sets in the original result.

 \begin{theorem}\label{res:eleven_sets_equal}
 Let $E$ be a vector lattice with the countable sup property, and suppose that $\ocdual{E}$ separates the points of $E$. Then $E$ admits a \uppars{necessarily unique} Hausdorff \uoLtop\ $\uoLt_E$. Take an ideal $I$ of $\ocdual{E}$ that separates the points of $E$, and take a vector sublattice $F$ of $E$. Then the following eleven vector sublattices of $E$ are all equal:
 \begin{enumerate}
 	\item $\sadh{\uoLt_E}{F}$ and $\overline{F}^{\sigma\uoLt_E}$;
 	\item $\suoadh{F}$ and $\suoclos{F}$;
 	\item $\uoadh{F}$ and $\uoclos{F}$;
 	\item $\overline{F}^{\uoLt_E}$, $\overline{F}^{\abs{\sigma}(E,I)}$, and $\overline{F}^{\sigma(E,I)}$;
 	\item\label{part:double_order_adherence_and_order_closure} $(\oadh{\oadh{F}})$ and $\oclos{F}$.
 \end{enumerate}
 \end{theorem}

The equality of $\uoadh{F}$, $\oadh{\oadh{F}}$, and $\overline{F}^{\sigma(E,I)}$ can already be found in \cite[Theorem~2.2]{gao_leung:2018}, where it also noted that these sets coincide with the smallest order closed vector sublattice of $E$ containing $F$.

 \begin{proof}
 The equality of the first seven subsets is clear from \cref{res:seven_sets_equal}. Since we know from \cref{res:uo-lebesgue_topology_originating_from_separating_order_continuous_dual_of_order_dense_ideal} that $\uoLt_E=\unb_E\abs{\sigma}(E,I)$, it follows from \cite[Proposition~2.12]{taylor:2019} that   $\overline{F}^{\uoLt_E}=\overline{F}^{\abs{\sigma}(E,I)}$. Furthermore, from Kaplan's theorem (see \cite[Theorem~2.33]{aliprantis_burkinshaw_LOCALLY_SOLID_RIESZ_SPACES_WITH_APPLICATIONS_TO_ECONOMICS_SECOND_EDITION:2003}, for example) we know that $E$, when supplied with the Hausdorff locally convex $\abs{\sigma}(E,I)$-topology, has the same topological dual as when it is supplied with the Hausdorff locally convex $\sigma(E,I)$-topology. By the convexity of $F$, we have $\overline{F}^{\abs{\sigma}(E,I)}=\overline{F}^{\sigma(E,I)}$. This argument was already used in \cite[Proof of Lemma~2.1]{gao_leung:2018}.

 We turn to the two sets in part~\ref{part:double_order_adherence_and_order_closure}. It was established in \cite[Lemma~2.1]{gao_leung:2018} that $
 \uoadh{F}\subseteq\oadh{\oadh{F}}$; this is, in fact, valid for vector sublattices of general vector lattices. It was also observed there that, obviously, the fact that $I\subseteq\ocdual{E}$ implies that $\overline{F}^{\sigma(E,I)}$ is o-closed. Using also that we already know that $\uoadh{F}=\uoclos{F}$, we therefore have the following chain of inclusions:
 \[
 \uoclos{F}=\uoadh{F}\subseteq \oadh{\oadh{F}}\subseteq \oclos{F}\subseteq \overline{F}^{\sigma(E,I)}.
\]
 Since we also already know that $\uoclos{F}=\overline{F}^{\sigma(E,I)}$, the proof is complete.
 \end{proof}

\subsection*{Acknowledgements} During this research, the first author was supported by a grant of China Scholarship Council (CSC). The authors thank Jurie Conradie, Mitchell Taylor, and Vladimir Troitsky for helpful discussions. They are grateful to the anonymous referees for their careful reading of the manuscripts and the comments they made. These have led, amongst others, to the inclusion of a reference to \cite[Theorem~6.45]{taylor_THESIS:2018} and to \cref{res:order_separability_equivalences} now being established for vector sublattices rather than just for ideals as in the original manuscript.

\bibliographystyle{plain}
\urlstyle{same}
\bibliography{../../../../../../tex/templates/bibliography/general_bibliography}

\end{document}